\documentclass[11pt,a4paper,twoside]{amsart}
\usepackage{amssymb}
\usepackage[mathscr]{euscript}

\newtheorem*{Theoreme}{Th\'eor\`eme}
\newtheorem*{Corollaire}{Corollaire}  
\newtheorem*{Proposition}{Proposition}
\newtheorem*{Lemme}{Lemme}

\theoremstyle{definition}

\theoremstyle{remark}
\newtheorem*{Remarque}{Remarque}

\def\goth{\mathfrak }
\def\scr{\mathcal }
\def\Bbb{\mathbb }

\def\tl{\widetilde }

\def\ol{\overline }
\def\fl{\rightarrow }
\def\ffl{\mapsto }
\def\ssi{\Leftrightarrow }
\def\impl{\Rightarrow }
\def\x{\times }
\def\[{[\! [}
\def\]{]\!]}

\def\End{\text{\rm End}}
\def\GL{\text{\rm GL}}

\def\Gal{{\scr Gal }}

\def\id{\text{\rm id}}
\def\Ker{\text{\rm Ker }}
\def\im{\text{\rm Im }}

\def\Det{\text{\rm d\'et}}
\def\tr{\text{\rm tr}\, }

\def\diag{\text{\rm diag}}

\def\val{\,\text{\rm val}\, }

\def\ind{\text{\rm Ind}}
\def\indu{\text{\rm ind}}

\def\Cl{{\scr Cl} }
\def\Nc{{\scr N}_\tau }
\def\N{N_\tau}

\def\J{\Bbb J}
\def\K{\Bbb K}

\numberwithin{equation}{subsection}

\begin{document}

\title[Changements de base de $U(1,1)(F_{0})$]{\sc Changements de base explicites des repr\'esentations supercuspidales de $U(1,1)(F_0)$}
 
\author{Laure Blasco}

\address{D\'epartement de Math\'ematiques et U.M.R. 8628 du C.N.R.S.,
Universit\'e Paris-Sud, B\^atiment 425, 91405 Orsay cedex, France}

\email{Laure.Blasco@math.u-psud.fr}

\subjclass[2000]{Primary: 22E50, Secondary: 11F70}

\date{7 septembre 2009}
\maketitle
\begin{abstract} Soit $F_0$ un corps local non archim\'edien de caract\'eristique nulle et de ca\-rac\-t\'eristique r\'esiduelle impaire.  On d\'ecrit explicitement les changements de base des repr\'esentations supercuspidales de $U(1,1)(F_{0})$. C'est une \'etape vers la description du changement de base des paquets endoscopiques supercuspidaux de $U(2,1)(F_{0})$.
\end{abstract}

Bien que le titre n'en dise rien, cet article fait suite \`a \cite{Bl2} dans lequel nous d\'ecrivons explicitement le changement de base stable de certaines repr\'esentations supercuspidales du groupe unitaire $U(2,1)(F_{0})$ relativement \`a $F$ au groupe lin\'eaire $GL(3,F)$ o\`u $F_{0}$ est un corps $p$-adique de caract\'eristi\-que r\'esiduelle impaire et $F$ une extension quadratique de $F_{0}$. Cette description repose sur la classification des repr\'esentations supercuspidales par la th\'eorie des types de C. Bushnell et Ph. Kutzko et sur les travaux de J. Rogawski \cite{Ro}. 

\smallskip	

Plus pr\'ecis\'ement, dans \cite{Bl2}, n'est explicit\'e que le changement de base stable de paquets cuspidaux de cardinal 1 (et on pense qu'ils y sont tous). Il reste donc \`a faire l'analogue pour les paquets cuspidaux endoscopiques, ce qui se r\'ealise en deux temps \cite[Ch.4, \S 2]{Ro} :
\begin{enumerate}
\item d\'eterminer les images par l'application de transfert des repr\'esenta\-tions de carr\'e int\'egrable du groupe $U(1,1)(F_{0})\x U(1)(F_{0})$ au groupe $U(2,1)(F_{0})$ ;
\item d\'ecrire le changement de base ``labile'' (ou ``instable'') des repr\'esen\-ta\-tions supercuspidales de $U(1,1)(F_{0})$ au groupe $GL(2,F)$.
\end{enumerate}
On conclut alors gr\^ace aux r\'esultats de J. Rogawski \cite[prop. 13.2.2(c)]{Ro}. Notons que pour la repr\'esentation de Steinberg de $U(1,1)(F_{0})$ et  ses tordues par un caract\`ere, la deuxi\`eme \'etape est faite dans \cite{Ro} (prop. 11.4.1 et d\'e\-monstration de la prop. 12.4.1).

\smallskip	

Ce texte pr\'esente la r\'ealisation du point (2). Elle repose sur la liste des types de $U(1,1)(F_{0})$ obtenue dans l'annexe de \cite{Bl1}, sur les r\'esultats concernant $GL(2,F)$ expos\'es dans les chapitres 5 et 8 de \cite{BH3} et sur ceux du chapitre 11 de \cite{Ro}. Comme dans le cas de $U(2,1)(F_{0})$, on d\'etermine le changement de base stable dont l'identit\'e de caract\`eres ne d\'epend que des classes de conjugaison et conjugaison tordue stables, le changement de base labile s'en d\'eduisant ais\'ement puisqu'il est explicitement reli\'e au pr\'ec\'edent (\cite[\S 11.4]{Ro} ou voir \S \ref{notations}). La compatibilit\'e du changement de base \`a la torsion par un caract\`ere permet de se restreindre aux repr\'esentations supercuspidales de niveau minimal parmi leurs tordues par un caract\`ere. 

\smallskip	

Dans le cas de $U(1,1)(F_{0})$, les paquets supercuspidaux sont bien connus~:  on en donne une description ``typique'' dans le paragraphe \ref{paquets}. Les paquets endoscopiques sont pr\'ecis\'ement les paquets de cardinal 2 \cite[prop. 11.1.1(a)]{Ro} et sont form\'es, en niveau strictement positif, de repr\'esentations {\it cuspidales scind\'ees} (c'est-\`a-dire celles dont le type provient d'une strate gauche fondamentale scind\'ee suivant un sous-groupe de Levi non rationnel sur $F_{0}$). Les repr\'esentations supercuspidales de niveau 0 sont toutes dans des paquets de cardinal 2 si $F$ est non ramifi\'ee sur $F_{0}$ et toutes sauf deux dans des paquets singletons si $F$ est ramifi\'ee sur $F_{0}$.  
 
Par contre, il est plus difficile de distinguer les ima\-ges des deux changements de base, toutes deux form\'ees de repr\'esentations admissibles, inva\-rian\-tes sous l'action du groupe de Galois de  $F/F_{0}$ et de caract\`ere central trivial sur $F_{0}^\x$. Pour un paquet endoscopique, emprunter la ``voie'' par l'application transfert de $U(1)(F_{0})\x U(1)(F_{0})$ \`a $U(1,1)(F_{0})$ et utiliser le changement de base de $U(1)(F_{0})$ \`a $F^\x$ permet d'ignorer cette question. Mais pour un paquet singleton o\`u l'on \'etablit le changement de base en v\'erifiant une identit\'e de caract\`eres entre repr\'esentations de $U(1,1)(F_{0})$ et de $GL(2,F)$ (voir (\ref{id})), rep\'erer les repr\'esentations supercuspidales de $GL(2,F)$ qui appartiennent \`a l'image du changement de base stable est un point crucial. 

Dans le paragraphe \ref{paquetsH}, on d\'etermine le changement de base stable des paquets endoscopiques en empruntant la m\'ethode de Y. Flicker \cite{Fl} puis on conclut gr\^ace \`a \cite[prop. 11.4.1(a)]{Ro}. 
Dans le paragraphe \ref{paquetssing}, on d\'ecrit les changements de base des paquets non endoscopiques en suivant la m\^eme d\'emarche que celle expos\'ee dans \cite{Bl2} dont on ne reprend pas tous les d\'etails. 
L'ensemble des r\'esultats est pr\'esent\'e aux corollaire \ref{BC-endoscopiques} et th\'eor\`eme \ref{BCdim2}.

Le dernier paragraphe expose un calcul technique permettant de comparer les caract\`eres de repr\'esentations d'un groupe et de l'un de ses sous-groupes. Il est utilis\'e \`a plusieurs reprises dans les paragraphes pr\'ec\'edents.

\smallskip	

Les m\'ethodes utilis\'ees n'ont pas d'originalit\'e mais fournissent les r\'esultats d\'esir\'es, non \'ecrits jusque-l\`a, sous une forme  que nous pourrons exploiter ult\'erieurement et sous un minimum d'hypoth\`eses ($F_{0}$ de caract\'eristique nulle et de caract\'eristique r\'esiduelle impaire). Pour ces raisons, on inclut les repr\'esenta\-tions cuspidales de niveau 0 d\'ej\`a \'etudi\'ees par J. Adler et J. Lansky (\cite{AL} et \cite{AL2}). 

\smallskip	
Je remercie Corinne Blondel pour ses encouragements patients et sa relecture critique.

\section{Notations}\label{notations}
Soient $F_0$ un corps local non archim\'edien, de caract\'eristique nulle 
et de caract\'eristique r\'esiduelle diff\'erente de deux et $F$ une extension
quadratique (s\'eparable) de $F_0$ d'indice de ramification $e_0$ et dont le
groupe de Galois est not\'e  $\Gamma$~: $\Gamma=\{1, {}^-  \}$.

On d\'esigne par ${\goth o}_0$ (resp. ${\goth o}$) l'anneau des entiers de
$F_0$ (resp. $F$), ${\goth p}_0$ (resp. ${\goth p}$) l'id\'eal maximal de
${\goth o}_0$ (resp. ${\goth o}$) et $\varpi _0$ (resp. $\varpi$) une
uniformisante de ${\goth p}_0$  (resp. ${\goth p}$). On choisit les
uniformisantes $\varpi$ et $\varpi _0$ telles que : $\varpi =\varpi _0$ si
$e_0=1$, $\varpi $ est de trace nulle et de norme $\varpi _0$ si $e_0=2$. On
note $k_0$ et $k$ les corps r\'esiduels de $F_0$ et $F$ respectivement et
$q$ le cardinal de $k_0$. 

Pour une extension $E$ de $F_0$, on conserve les m\^emes notations que
pour $F$, cette fois index\'ees par $E$. Si $L$ est une  sous-extension de $E$, on d\'esigne par $E_{\vert L}^1$ le groupe des \'el\'ements de $E$ dont la norme dans $L$ est 1.

\smallskip

On fixe un caract\`ere additif $\psi_0$ de $F_0$, de conducteur ${\goth
p}_{0}$. Sa compos\'ee avec la trace $\tr _{F/F_0}$ est un caract\`ere
$\psi$ de $F$ de conducteur ${\goth p}$. On fixe \'egalement un prolongement \`a $F^\times$ du caract\`ere $\omega_{F/F_0}$ de $F_0^\times$ associ\'e \`a $F/ F_0$ par la th\'eorie du corps de classes. On choisit ce prolongement, not\'e $\mu$,
 \'egal \`a $x\ffl (-1)^{\val_F(x)}$ si $F$ n'est pas ramifi\'ee sur $F_0$ ;  trivial sur $1+{\goth p}$ et \'egal \`a la constante de Langlands $\lambda_{F_{\vert F_{0}}}(\psi_{0})$ en  $\varpi$ si
$F$ est ramifi\'ee sur $F_0$ (voir par exemple \cite[\S 34.3]{BH3}.

Si $\chi$ est un caract\`ere de $F_{\vert F_{0}}^1$, on d\'esigne par $\tl \chi$ le caract\`ere de $F^\x$ d\'efini par~: $\tl \chi (x)=\chi(\frac{x}{\ol x}), x\in F^\x$. L'application $\chi \mapsto \tl \chi$ n'est autre que le changement de base stable de $U(1)(F_{0})$ \`a $GL(1,F)$.

\medskip

Soient $V$ un $F$-espace vectoriel de dimension 2 muni d'une forme
hermitien\-ne ``isotrope'' non d\'eg\'en\'er\'ee $< , >$, $G=U(1,1)(F_0)$ son groupe
d'isom\'etries et $\widetilde G=GL (2,F)$ son groupe d'automorphismes. On note $Z$ le centre de $U(1,1)$ :  $Z(F_{0})$ s'identifie \`a $F_{\vert F_{0}}^1$ et $Z(F)$ \`a $F^\x$.\\
Le groupe $\Gamma$ agit sur $\widetilde G$ par : l'\'el\'ement
non trivial de $\Gamma$ transforme un \'el\'ement $g$ de $\widetilde G$ en $\tau
(g):= \sigma(g)^{-1}$ o\`u $\sigma$ d\'esigne l'involution d\'efinie sur $\End
_FV$ et associ\'ee \`a $<,>$. Alors, $G$ n'est autre que le groupe des points
de $\widetilde G$ fixes sous $\Gamma$ : $G=\widetilde G^\tau$.

\smallskip	

\noindent On fixe une base hyperbolique ${\scr B}=(e_{-1}, e_{1})$ de $V$ et on note $g_{0}$ la similitude de matrice $\left(\begin{array}{cc} 1&0\\0&\alpha_{0}\end{array}\right)$ o\`u $\alpha_{0}=\varpi_{0}$ si $F$ n'est pas ramifi\'ee, $\alpha_{0}\in {\goth o}_{0}^\x$ n'est pas une norme si $F$ est ramifi\'ee. Alors le groupe  $GU(1,1)(F_{0})$ des similitudes unitaires de $V$ est la r\'eunion de $Z(F) G$ et $g_{0}Z(F) G$.

\smallskip

\noindent On fixe \'egalement un \'el\'ement $\varepsilon$ de $F$ : $\varepsilon$ est une unit\'e de trace nulle si $F$ n'est pas ramifi\'ee sur $F_{0}$, $\varpi$ sinon. On note \'egalement $\varepsilon$ l'\'el\'ement de $\tl G$ dont la matrice dans ${\scr B}$ est $\left( \begin{array}{cc} 1&0\\0&\varepsilon \end{array}\right)$. Alors l'application $\Phi$ d\'efinie par :
$$g \mapsto {}^\varepsilon g:=\left( \begin{array}{cc} 1&0\\0&\varepsilon \end{array}\right) g \left( \begin{array}{cc} 1&0\\0&\varepsilon \end{array}\right)^{-1}$$
 est un isomorphisme  entre $GZ(F)$ et $GL(2,F_{0})^+Z(F)$ qui identifie $SU(1,1)(F_{0})$ et $SL_{2}(F_{0})$. Rappelons que : $GL(2,F_{0})^+=\{ g\in GL(2,F_{0}) \vert \Det g\in N_{F\vert F_{0}}(F^\x)\}$.
 
 \smallskip	
 
\noindent Le groupe endoscopique elliptique $H=U(1)(F_{0})\x U(1)(F_{0})$ de $G$ s'identifie au sous-groupe  des \'el\'ements de $G$ dont la matrice dans la base orthogonale ${\scr B}^o=(e_{-1}+\frac{1}{2}e_{1},e_{-1}-\frac{1}{2}e_{1})$ est diagonale. A conjugaison pr\`es dans $G$, il existe au plus deux plon\-gements de $H$. Ils sont alors conjugu\'es par $g_{0}$.

\medskip	

Les d\'efinitions de l'application de transfert de $H$ \`a $G$ (cf. prop. \ref{donnees}) et du changement de base labile \cite[\S 4.7] {Ro} d\'ependent du choix d'un  caract\`ere de $F^\x$ prolongeant $\omega_{F/F_{0}}$.  On choisit le caract\`ere $\mu$ d\'efini pr\'ec\'edemment. On a alors \cite[\S 11.4]{Ro} :
{\it si $\pi$ est une repr\'esentation admissible de $G$, son image par le changement de base stable est
$\tl \pi$ si et seulement si son image par le changement de base labile  est $\tl \pi \cdot \mu\circ \Det$.}

\smallskip	

Pour \'etablir le changement de base stable des paquets supercuspidaux singletons de $G$, on v\'erifie l'identit\'e de caract\`eres qui le caract\'erise et qui s'exprime gr\^ace \`a la  {\it  norme cyclique}, not\'ee  $\Nc $ (\ref{id}). Il s'agit d'une bijection de l'ensemble des classes de $\tau$-conjugaison stable de $\tl G$ dans l'ensemble des classes de conjugaison stable de $G$, qui associe \`a la classe de $\tau$-conjugaison stable de $g$, l'intersection de la classe de $\tl G$-conjugaison de $\N (g)=g\tau(g)$ avec $ G$ \cite{Ko}.\\
Rappelons que la classe de $\tau$-conjugaison d'un  \'el\'ement $g$ de $\tl G$, ${\scr Cl}_\tau (g)$, est l'ensemble des \'el\'ements de $\tl G$ de la forme $h^{-1}g\tau(h)$ avec $ h\in {\tl G}$ et que sa classe de $\tau$-conjugaison stable $\Cl_{\tau}^{st} (g)$ est l'ensemble des \'el\'ements $g'$ de $\tl G$ tel que $\N (g')$ soit $\tl G$-conjugu\'e \`a $\N (g)$.\\
 De fa\c con analogue, pour un \'el\'ement $x$ de $G$, on d\'esigne par $\Cl (x)$ sa classe de conjugaison et par $\Cl^{st}(x)$ sa classe de conjugaison stable (c'est-\`a-dire de conjugaison sous $\tl G$).

\section{Les paquets supercuspidaux de $U(1,1)(F_{0})$.}\label{paquets}  Le groupe $GU(1,1)(F_{0})$ agit par conjugaison sur l'ensemble des repr\'esen\-ta\-tions lisses irr\'educ\-ti\-bles de $U(1,1)(F_{0})$. Les orbites sont exac\-tement les paquets de $U(1,1)(F_{0})$ \cite[\S 11.1]{Ro}. On en d\'eduit : 

\begin{Proposition}\mbox{}

\noindent (i)  Les repr\'esentations irr\'eductibles tr\`es cuspidales de $U(1,1)(F_{0})$ de niveau strictement positif sont seules dans leur paquet. Les autres repr\'esentations cuspidales de niveau strictement positif appartiennent \`a des paquets de cardinal $2$.\\(ii) Si $F$ est non ramifi\'ee sur $F_{0}$, tous les paquets form\'es par des repr\'esen\-tations irr\'eductibles cuspidales de niveau $0$ sont de cardinal 2.\\
(iii) Si $F$ est ramifi\'ee sur $F_{0}$, tous les paquets form\'es par des repr\'esenta\-tions irr\'eductibles cuspidales de niveau $0$ sont de cardinal $1$ sauf un. Ce dernier est $\{ \ind_{K}^G\sigma^+,\ind_{K}^G\sigma^-\}$ o\`u $K$ est le sous-groupe parahorique maximal et $\sigma^\pm$ les deux repr\'esentations de $SL_{2}(k_{0})$ de dimension $\frac{q-1}{2}$ \cite{Sp}. 
\end{Proposition}

La suite du paragraphe  justifie ces assertions.

\smallskip	

Soient $\pi$ une repr\'esentation irr\'eductible cuspidale de $G$ et  $(J,\lambda)$ un type pour cette repr\'esentation : $\pi=\ind_{J}^G\lambda$. Il suffit d'\'etudier si $\pi$ et $^{g_{0}}\pi$, sa conjugu\'ee par $g_{0}$, sont isomorphes. On distingue deux cas selon le niveau.	

\subsection{$\pi$ est de niveau strictement positif.} Alors $(J,\lambda)$ provient d'une strate gauche semi-simple ${\goth s}=(\scr L, n,\left[\frac{n}{2}\right],b)$  (au sens de \cite[d\'ef. 4.9]{St}). On note $Z_{\tl G}(b)$ le centralisateur de $b$ dans $\tl G$ et $Z_{G}(b)$ celui de $b$ dans $G$ :  $Z_{G}(b)=Z_{\tl G}(b)\cap G$.

\smallskip	

Si $\pi$ est tr\`es cuspidale \cite[\S A.5]{Bl1}, $Z_{\tl G} (b)\cup \{0\}$ est l'extension quadratique $F[b]$ de $F$. L'\'el\'ement $c=b-\frac{1}{2}\tr_{F[b]_{\vert F}} b$  engendre aussi $F[b]$ sur $F$, appartient \`a l'alg\`ebre de Lie de $G$ et engendre une extension quadratique $F_{0}[c]$ sur $F_{0}$, autre que $F$. On v\'erifie que les \'el\'ements de $F_{0}[c]^\x$ sont des \'el\'ements de $GU(1,1)(F_{0})$ dont le rapport est \'egal \`a leur norme dans $F_{0}$. Il en existe donc un, $c_{0}$, qui appartient \`a $g_{0}F_{0}^\x U(1,1)(F_{0})$. Alors : $^{g_{0}}\pi\simeq {}^{c_{0}}\pi\simeq \pi$.

\smallskip	

Si $\pi$ est cuspidale scind\'ee \cite[\S A.6]{Bl1}, $Z_{\tl G}(b)$ est isomorphe \`a $F^\x \x F^\x$. On montre par l'absurde que l'entrelacement ${\scr E}({\goth s},{}^{g_{0}}{\goth s})$ des strates semi-simples ${\goth s}$ et ${}^{g_{0}}{\goth s}=(g_{0 }\scr L, n,\left[\frac{n}{2}\right],g_{0}^{-1}bg_{0})$ contenues dans $\lambda$ et ${}^{g_{0}}\lambda$ respectivement,  est vide (d'o\`u l'on d\'eduit que $\lambda$ et ${}^{g_{0}}\lambda$ ne sont pas entrelac\'ees dans $G$).\\ 
L'entrelacement  $\tl {\scr E}({\goth s},{}^{g_{0}}{\goth s})$ de ${\goth s}$ et ${}^{g_{0}}{\goth s}$ dans $\tl G$  est \'egal \`a $\tl {\scr E}({\goth s})g_{0}^{-1}$ o\`u $\tl {\scr E}({\goth s})$ est l'entrelacement de ${\goth s}$ dans $\tl G$. De plus, par le th\'eor\`eme 4.10 pr\'ecis\'e par  (4.15)  de \cite{St}, $\tl {\scr E}({\goth s})$ est \'egal \`a $\tl U_{m}(\scr L)Z_{\tl G}(b) \tl U_{m}(\scr L)$ avec $m=\left[\frac{n+1}{2}\right]$.
On d\'eduit que~:
\begin{equation*}\begin{split}
{\scr E}({\goth s},{}^{g_{0}}{\goth s})&=\tl {\scr E}({\goth s},{}^{g_{0}}{\goth s})^\tau=\left(\tl U_{m}(\scr L) Z_{\tl G}(b)g_{0}^{-1}\tl U_{m}(g_{0} {\scr L})\right)^\tau \\
&=\cup_{z\in Z_{\tl G}(b)}\left(\tl U_{m}(\scr L)  zg_{0}^{-1}\tl U_{m}(g_{0} {\scr L})\right)^\tau.
\end{split}\end{equation*}
Si ${\scr E}({\goth s},{}^{g_{0}}{\goth s}) \neq \emptyset$, il existe $z\in Z_{\tl G}(b)$ tel que : $\left(\tl U_{m}(\scr L)  zg_{0}^{-1}\tl U_{m}(g_{0} {\scr L})\right)^\tau\neq \emptyset$. Fixons un tel $z$. Comme $\scr L$ et $g_{0}\scr L$ sont des cha\^\i nes autoduales ($g_{0}$ appartient \`a $GU(1,1)(F_{0})$), les deux pro-$p$-sous-groupes $\tl U_{m}(\scr L)$ et $\tl U_{m}(g_{0}\scr L)$ sont invariants par $\tau$. Par le lemme 2.2 de \cite{St}, g\'en\'eralis\'e sans encombre au cas de deux pro-$p$-groupes, $\tl U_{m}(\scr L)  zg_{0}^{-1}\tl U_{m}(g_{0} {\scr L})$ est invariant par $\tau$. En particu\-lier $\tau(zg_{0}^{-1})$, qui est \'egal \`a $\alpha_{0}\tau(z)g_{0}^{-1}$, appartient \`a $\tl U_{m}(\scr L)  zg_{0}^{-1}\tl U_{m}(g_{0} {\scr L})=\tl U_{m}(\scr L)  z\tl U_{m}({\scr L})g_{0}^{-1}$ d'o\`u : $\alpha_{0}z^{-1}\tau(z)= z^{-1}uz\cdot u'$ avec $u, u'\in \tl U_{m}(\scr L)$.\\
Exprim\'es dans une base orthogonale form\'ee de vecteurs propres de $b$ et adapt\'ee \`a $\scr L$, le terme de gauche est une matrice diagonale dont les termes diagonaux appartiennent \`a $\alpha_{0}N_{F/F_{0}}(F^\x)$ tandis que le terme de droite a au moins un de ses coefficients diagonaux dans $1+{\goth p}^m$. L'\'egalit\'e entra\^\i ne donc que $\alpha_{0}$ doit \^etre une norme de $F^\x$ dans $F_{0}$ ce qui est absurde.

\subsection{$\pi$ est de niveau $0$.} Il existe un sous-groupe parahorique $K$ maximal (au sens de Bruhat-Tits) et une repr\'esentation irr\'eductible $\sigma$ de $K$, triviale sur le radical pro-unipotent $K_{1}$ et d\'efinissant une repr\'esentation cuspidale  du quotient r\'eductif (connexe) $K/K_{1}$, tels que $\pi$ soit isomorphe \`a $\ind_{K}^G\sigma$ (dans ce cas, $K$ est son propre normalisateur) \cite[\S\S 2 et 3]{Mo}.

\smallskip	

Quand $F$ est non ramifi\'ee sur $F_{0}$, il existe, \`a conjugaison pr\`es, deux sous-groupes parahoriques maximaux : $K=U_{0}(\scr L)$ o\`u $\scr L$ est une cha\^\i ne de r\'eseaux autoduale de p\'eriode 1, d'invariant pair ou impair, et $K'$ son conjugu\'e par $g_{0}$. Dans ce cas,  $\pi=\ind_{K}^G\sigma$ et ${}^{g_{0}}\pi=\ind_{K'}^G{}^{g_{0}}\sigma$  ne sont pas isomorphes. \\
En effet, puisque $U(1,1)(F_{0})=K SU(1,1)(F_{0})$, les restrictions \`a $SU(1,1)(F_{0})$ de $\pi$ et ${}^{g_{0}}\pi$ sont deux repr\'esentations irr\'educ\-ti\-bles  de $SU(1,1)(F_{0})$ dont les  conjugu\'ees par $\varepsilon$,
not\'ees ${}^\varepsilon \pi$ et ${}^\varepsilon ({}^{g_{0}}\pi)$ respectivement, sont deux repr\'esenta\-tions de $SL_{2}( F_{0})$, cuspidales, irr\'eductibles et de niveau 0 : ${}^\varepsilon \pi=\ind_{{}^\varepsilon K}^{SL_{2}(F_{0})} {}^\varepsilon \sigma$ et ${}^\varepsilon ({}^{g_{0}}\pi)={}^{g_{0}}({}^\varepsilon \pi)$ (puisque $\varepsilon $ et $g_{0}$ commutent).
Elles sont donc toutes deux sous-repr\'esentations de la restriction \`a $SL_{2}(F_{0})$ d'une  repr\'esentation irr\'eductible cuspidale de niveau 0 de $GL_{2}(F_{0})$ et ne sont donc pas \'equivalentes \cite[th. 4.4]{KS}.

\smallskip	

Lorsque $F$ est ramifi\'ee, il existe une unique classe de conjugaison de sous-groupes parahoriques maximaux repr\'esent\'ee par $U_{0}(\scr L)$ o\`u $\scr L$ est la  cha\^\i ne de r\'eseaux autoduale de p\'eriode 1 et d'invariant impair. Le quotient $U_{0}(\scr L)/U_{1}(\scr L)$
est isomorphe \`a $SL_{2}(k_{0})$.\\
Ainsi, $\sigma$ d\'efinit une repr\'esentation cuspidale irr\'eductible de $SL_{2}(k_{0})$. D'apr\`es T. A. Springer \cite[II \S 3]{Sp}, $\sigma$ est de dimension $q-1$ ou $\frac{q-1}{2}$ et est d\'ecrite par son caract\`ere. Il est alors ais\'e de montrer que si $\sigma$ est de dimension $q-1$, elle est isomorphe \`a sa conjugu\'ee par $g_{0}$ tandis que si sa dimension est  $\frac{q-1}{2}$, sa conjugu\'ee est l'autre repr\'esentation cuspidale de m\^eme dimension, not\'ee $\sigma'$. Dans le premier cas, $\pi$ est isomorphe \`a sa conjugu\'ee par $g_{0}$. Dans le deuxi\`eme cas, montrons que $\sigma$ et $\sigma'$ ne sont pas entrelac\'ees.\\ 
Gr\^ace \`a la d\'ecomposition de Cartan, $K$ \'etant un bon compact dans ce cas, il suffit de montrer que $g=\left(\begin{array}{cc}\varpi&0\\ 0&{\ol \varpi} ^{-1}\end{array}\right)^m$, $m\in {\Bbb N}$, n'entrelace pas $\sigma$ et $\sigma'$. Si $m=0$, c'est clair car $\sigma$ et $\sigma'$ ne sont pas isomorphes. Si $m>0$, $K\cap {}^gK$ contient le sous-groupe $\left(\begin{array}{cc}1&0\\ \varpi{\goth o}_{0}&1\end{array}\right)$. Or la restriction de $\sigma'$ \`a ce sous-groupe ne contient pas le caract\`ere trivial tandis que celle de ${}^g\sigma$ est triviale.

\section{Les paquets endoscopiques cuspidaux et les caract\`eres de $H$.}\label{paquetsH}
Il s'agit d'iden\-tifier le caract\`ere $\theta$ de $H$ associ\'e \`a un paquet endoscopique cuspidal de $G$ par l'application de transfert. 

\smallskip	

Remarquons que lorsque  $F$ est ramifi\'ee sur $F_{0}$, il existe un unique paquet endoscopique de niveau 0 et deux caract\`eres r\'eguliers de $H$ de niveau 0 : $\theta= 1\otimes \chi$  et $\chi\otimes 1={}^\tau \theta$ o\`u $\chi$ est le caract\`ere de $F_{\vert F_{0}}^1$ d'ordre 2 trivial sur $(1+{\goth p})_{\vert F_{0}}^1$. N\'ecessairement, le paquet endoscopique est l'image de $\theta$ (et ${}^\tau \theta$) par l'application de transfert. Par la suite, on suppose donc :
$$ n>0 \quad \text{ ou } \quad F \text{ n'est pas ramifi\'ee sur } F_{0}.$$ 

\subsection{}\label{donnees} Soit $\Pi$ un paquet endoscopique cuspidal de $G$. Il est enti\`erement d\'etermin\'e par la donn\'ee d'un type simple $(J,\lambda)$ :
\begin{enumerate}
\item[(a)]   $J= HU_{m}(\scr L)$ o\`u $\scr L$ est la cha\^\i ne de ${\goth o}$-r\'eseaux autoduale stable par $H$, c'est-\`a-dire celle de p\'eriode 1, d'invariant pair, et $m=[\frac{n+1}{2}]$, $n\in \Bbb N$ ; 
\item[(b)] {\it Cas de niveau 0, $F$ non ramifi\'ee} : $\lambda$ est le rel\`evement \`a $U_{0}(\scr L)$ d'une repr\'esentation irr\'eductible cuspidale $\sigma$ de $U(1,1)(k_{0})$ v\'erifiant :
\item[]  $\quad \quad \forall x\in H(k_{0}), \tr \sigma(x)=\left\{ \begin{array}{ll} (q-1) \cdot \chi(x) & \text{ si } x\in Z(k_{0})\\ -(\chi (x)+\chi(wxw^{-1}))& \text{ si } x\not \in Z(k_{0})\\ \end{array}\right.$
\item[] pour un caract\`ere r\'egulier $\chi$ de $H(k_{0})$ \cite[\S 6]{En}. L'\'el\'ement $w$ appartient au normalisateur de $H(k_{0})$ mais non \`a $H(k_{0})$. Notons $\theta_{1}\otimes \theta_{2}$ le caract\`ere de $H$ relevant $\chi$.
\item[(b')] {\it Cas de niveau strictement positif} : $\lambda$ est une repr\'esentation irr\'educti\-ble de $J$ dont la restriction \`a $U_{m}(\scr L)$ est multiple d'un caract\`ere $\psi_{b}$ avec $b\in   {\goth a}_{-n}(\scr L)^-\setminus  {\goth a}_{-n+1}(\scr L)^-$ et dont la restriction \`a $H$ v\'erifie : 
\item[]  $\quad \quad \forall x\in H, \tr \lambda(x)=\left\{ \begin{array}{ll} \dim \lambda \cdot \theta_{1}\otimes \theta_{2}(x) & \text{ si } x\in Z(F_{0})\\ -\theta_{1}\otimes \theta_{2}(x) & \text{ si } x\not \in Z(F_{0})\\ \end{array}\right.$
\item[] pour un caract\`ere r\'egulier $\theta_{1}\otimes \theta_{2}$ de $H$ prolongeant ${\psi_{b}}_{\vert H\cap U_{m}(\scr L)}$ (voir \S \ref{Traces} et en particulier  \S \ref{precisions}, (1)). Pr\'ecisons que $b$ est de la forme $b_{1}\oplus b_{2}$ dans la d\'ecomposition de $F^2$ d\'efinie par ${\scr B}^o$ avec $b_{1}-b_{2}\not \in {\goth p}^{-n+1}$ (quitte \`a tordre $\Pi$ par un caract\`ere de $G$).
\end{enumerate}
On a alors : $\Pi=\{ \pi, {}^{g_{0}}\pi\}$ o\`u $\pi=\ind_{J}^G \lambda$.

\begin{Proposition}[\cite{Ro}] Soit $\theta$ le caract\`ere de $H$ correspondant \`a $\Pi$ par l'application de transfert. Il existe un unique  $\epsilon_{\theta}\in \{\pm 1\}$ tel que pour tout $h\in H_{G-reg}$, $h=(h_{1}, h_{2})$ (c'est-\`a-dire $h$ de matrice $\diag(h_{1}, h_{2})$ dans ${\scr B}^o$),
\begin{equation}\label{IC-tr}
\tr{\pi}(h)-\tr{{}^{g_{0}}\pi}(h)=\epsilon_{\theta}\frac{\iota(h)}{D_{G}(h)}(\theta(h)+\theta(whw^{-1}))
\end{equation}
o\`u $\iota(h)=\mu(h_{1}-h_{2})$, $D_{G}(h)=\left\vert \frac{(h_{1}-h_{2})^2}{h_{1}h_{2}}\right\vert_{F_{0}}^{\frac{1}{2}}$ et $w=\left(\begin{array}{cc} 0&1\\1&0Ê\end{array}\right)$ dans ${\scr B}^o$.
\end{Proposition}

\subsection{} Evaluons le membre de gauche de (\ref{IC-tr}) en suivant la m\'ethode de Y. Flicker \cite{Fl}, c'est-\`a-dire en se ramenant aux r\'esultats de \cite[\S 2]{LL}, gr\^ace \`a l'isomorphisme $\Phi$ entre $GZ(F)$ et $GL(2,F_{0})^+Z(F)$ (cf. \S \ref{notations}). Pr\'ecisons que dor\'enavant, toutes les matrices sont exprim\'ees relativement \`a  la base $\scr B$. \\
On note $i$ le plongement de $F^\times$ dans $GL(2,F_{0})^+$ qui \`a $x=\alpha+\varepsilon \beta\in F^\x$ ($\alpha,\beta\in F_{0}$) associe $i(x)=\left( \begin{array}{cc} \alpha&2\beta\\ \frac{\varepsilon^2}{2}\beta&\alpha\end{array}\right)$.

\begin{Lemme}\label{Phi}\mbox{}\\
(i) Soient $h=(h_{1}, h_{2})\in H$, $x\in F^\x$ et $z\in Z(F)$. Alors :
\begin{align*}
&\Phi (hz)=i(\ol a h_{1})\cdot \ol a^{-1}z \text{ o\`u } a\in F^\x \text{ tel que } a\ol a^{-1}=\Det h=h_{1}h_{2}\\\text{ et }\quad &\Phi^{-1} (xz)=h_{x}\cdot \ol xz \text{ o\`u } h_{x}=(x\ol x^{-1},1)\in H.
\end{align*}
En particulier, $\Phi$ identifie $HZ(F)$ et $i(F^\x) Z(F)$.\\
(ii)  Soient $k\in {\Bbb N}$ et $\scr L^0$ la cha\^\i ne de ${\goth o}_{0}$-r\'eseaux dans $F_{0}^2$ stable par $F^\x$. Lorsque $k=0$, on suppose en outre que $F$ n'est pas ramifi\'ee sur $F_{0}$. Alors ~: \\
$\Phi(U_{k}({\scr L})Z(F_{0}))= \{ uz, u\in U_{k}({\scr L}^0), z\in Z(F) \text{ tels que } \Det u\cdot N_{F_{\vert F_{0}}}(z)=1\}$.\\
(iii) $\Phi( JZ(F))\cap GL(2,F_{0})^+=i(F^\x) U_{m}({\scr L}^0)$.
\end{Lemme}

\begin{proof} L'assertion {\it (i)} est obtenue par de simples calculs tandis que {\it (iii)} est une cons\'equence imm\'ediate de {\it (i)} et {\it (ii)}.\\
 {\it (ii)} Supposons d'abord $k>0$. 
Soient $x\in U_{k}(\scr L)$ et $z\in Z(F_{0})$. Il existe $y\in 1+{\goth p}^k$ tel que : $y\ol y^{-1}=\Det x$, puis $x'\in U_{k}(\scr L)\cap SU(1,1)(F_{0})$ tel que : $x=\left( \begin{array}{cc} y&0\\ 0&\ol y^{-1}\end{array}\right)x'$. Alors : $\Phi(xz)=\left( \begin{array}{cc} y\ol y&0\\ 0&1\end{array}\right)\Phi(x')\cdot \ol yz$.\\
Or $x'$ est de la forme $\left( \begin{array}{cc} a&\varepsilon b\\ \varepsilon^{-1}c& d\end{array}\right)$ avec $a,d\in (1+{\goth p}^k)\cap F_{0}=1+{\goth p}_{0}^{[\frac{k+1}{e_{0}}]}$ et $\varepsilon b, \varepsilon^{-1}c\in {\goth p}^k \cap F_{0}={\goth p}_{0}^{[\frac{k}{e_{0}}]}$ donc $\Phi(x')$ appartient \`a 
\begin{equation*}\left\{\begin{aligned}
\left( \begin{array}{cc} 1+{\goth p}_{0}^k&{\goth p}_{0}^k\\{\goth p}_{0}^k&1+{\goth p}_{0}^k \end{array}\right)\cap SL(2)(F_{0})\subset U_{k}({\scr L}^0) & \text{ si } e_{0}=1 \\ 
\left( \begin{array}{cc} 1+{\goth p}_{0}^{[\frac{k+1}{2}]}&{\goth p}_{0}^{[\frac{k}{2}]}\\{\goth p}_{0}^{[\frac{k}{2}]+1}&1+{\goth p}_{0}^{[\frac{k+1}{2}]} \end{array}\right)\cap SL(2)(F_{0})\subset U_{k}({\scr L}^0) & \text{ si } e_{0}=2.
\end{aligned}\right.
\end{equation*}
De plus, $y\ol y\in 1+\tr {\goth p}^k$ donc $\left( \begin{array}{cc} y\ol y&0\\ 0&1\end{array}\right)\in U_{k}({\scr L}^0)$. Ainsi, $\Phi(x)$ appartient \`a $U_{k}({\scr L}^0) Z(F)$ et $\Det\left(\left( \begin{array}{cc} y\ol y&0\\ 0&1\end{array}\right)\Phi(x')\right)\cdot N_{F_{\vert F_{0}}}(\ol y^{-1}z)=1$.\\
L'inclusion inverse se montre de fa\c con analogue. \\
Lorsque $k$ est nul, $F$ est non ramifi\'ee sur $F_{0}$ : les m\^emes arguments fournissent les m\^emes r\'esultats.
\end{proof}

\subsection{}
Pour poursuivre, on choisit un caract\`ere $\Omega$ de $Z(F)$ prolongeant le carac\-t\`ere central  $\omega_{\pi}$ de $\pi$ et tel que 
\begin{equation}\label{Omega}
\left\{ \begin{aligned} 
\forall u\in {\goth a}_{[\frac{n}{2}]+1}(\scr L)\cap Z(F), \quad \Omega(1+u)=\psi((b_{1}+b_{2})u) &\text{ si } n>0\\
\Omega=\Theta_{1}\Theta_{2}\quad \text{o\`u $\Theta_{i}$ est un prolongement de $\theta_{i}$ \`a $F^\x$} & \text{ si } n=0.\end{aligned}\right.
\end{equation}
On prolonge alors $\pi$ et ${}^{g_{0}}\pi$ \`a $GZ(F)$ en faisant agir $Z(F)$ via $\Omega$. On note encore $\pi$ et ${}^{g_{0}}\pi$ ces prolongements puis :
\begin{equation*}
\pi_{0}^+=\pi\circ \Phi^{-1}_{\vert GL_{2}(F_{0})^+}, \quad \pi_{0}^-={}^{g_{0}}\pi\circ \Phi^{-1}_{\vert GL_{2}(F_{0})^+}={}^{g_{0}}\pi_{0}^+ , \quad \pi_{0}=\ind_{GL_{2}(F_{0})^+}^{GL_{2}(F_{0})} \pi_{0}^+.
\end{equation*}

\begin{Lemme} Posons : $J_{0}=F^\x U_{m}({\scr L}^0)$, $J_{0,c}={\goth o}_{0}^\x U_{m}({\scr L}^0)$. Notons  $\Lambda$ la repr\'esentation $(\lambda\Omega)\circ \Phi^{-1}_{\vert J_{0}}$ de $J_{0}$ et $\Lambda_{c}$ sa restriction \`a $J_{0,c}$.\\
(i) La paire $(J_{0,c}, \Lambda_{c})$ est un type simple maximal de $\pi_{0}$ et $\pi_{0}=\ind_{J_{0}}^{GL_{2}(F_{0})} \Lambda$. La repr\'esentation $\pi_{0}$ est donc  irr\'eductible et cuspidale.\\
(ii) Cas $n>0$. La restriction de $\Lambda$ \`a $U_{[\frac{n}{2}]+1}(\scr L ^0)$ est un multiple du caract\`ere $\psi_{0, \alpha}$ o\`u $\alpha=i(b_{1}-b_{2})
\in {\goth a}_{-n}({\scr L}^0)\setminus {\goth a}_{-n+1}({\scr L}^0)$.\\
La restriction de $\Lambda$ \`a $F^\x$ v\'erifie pour tout $x\in F^\x$ :
\begin{itemize}
\item[] si $\dim \Lambda=1$, $\quad \Lambda (x)= \tl \theta_{1}(x)\Omega (\ol x)$ ;
\item[] si $\dim \Lambda>1$, $\quad \tr \Lambda (x)=\left\{ \begin{array}{ll} \dim \Lambda \cdot  \tl \theta_{1}(x)\Omega (\ol x)&\text{ si } x\in F_{0}^\x\\ - \tl \theta_{1}(x)\Omega (\ol x)&\text{ sinon.}\end{array}\right.$.
\end{itemize}
(iii) Cas $n=0$ ($F$ n'est pas ramifi\'ee). La restriction de $\Lambda$ \`a $U_{0}(\scr L ^0)$ est le rel\`evement  d'une repr\'esentation $\sigma_{0}$ de $GL_{2}(k_{0})$ caract\'eris\'ee par :
\begin{equation*}
\left\{
\begin{aligned}
\tr \sigma_{0}(z)=& (q-1)\Omega(z), \quad z\in k_{0}^\x\\
\tr \sigma_{0}(zn)=&-\Omega(z), \quad z\in k_{0}^\x, n\in N(k_{0})-\{ \id\} \\
\tr \sigma_{0}(x)=& -\Omega(\ol x)(\tl \theta_{1}(x)+\tl \theta_{2}(x)), \quad x\in k^\x-k_{0}^\x
\end{aligned}
\right.
\end{equation*}
et $\Lambda (\varpi_{0})=\Omega (\varpi_{0})$.
\end{Lemme}

\begin{proof}
La premi\`ere assertion est une cons\'equence des deux suivan\-tes et ces derni\`eres proviennent des propri\'et\'es de $\lambda$ et de $\Omega$ via $\Phi$.\\ 
Plus pr\'ecis\'ement, pour \'etablir {\it (ii)}, on \'etudie d'abord la restriction de $\Lambda$ \`a $U_{[\frac{n}{2}]+1}(\scr L ^0)\cap SU(1,1)(F_{0})$. Un calcul du m\^eme style que dans le lemme \ref{Phi} montre que cette restriction est multiple de $\psi_{0,\alpha}$ o\`u $\alpha=i(b_{1}-b_{2})+D$ avec $D$ une matrice diagonale. Ensuite, si $x=\diag(z,\ol z^{-1})\in U_{[\frac{n}{2}]+1}(\scr L ^0)$, $\Lambda (x)$ est multiple de $\psi_{{}^\varepsilon b}(x)\Omega(\ol z)$, qui vaut 1 par choix de $\Omega$ (\ref{Omega}). On peut donc prendre $D=0$. L'\'el\'ement $\alpha$ ainsi d\'efini appartient \`a ${\goth a}_{-n}({\scr L}^0)\setminus {\goth a}_{-n+1}({\scr L}^0)$ par \ref{donnees} (b'). La suite est imm\'ediate, tout comme l'assertion {\it (iii)}.
\end{proof}

\begin{Corollaire} Soit $\Delta_{\theta_{1}\otimes \theta_{2}}$ le caract\`ere de $F^\x$ de niveau $0$ d\'efini par : 
\begin{itemize}
\item[(i)] quand $F$ n'est pas ramifi\'ee sur  $F_{0}$, $\Delta_{\theta_{1}\otimes \theta_{2}}=\mu$ ;
\item[(ii)] quand $F$ est ramifi\'ee sur $F_{0}$, 

${\Delta_{\theta_{1}\otimes \theta_{2}}}_{\vert F_{0}^\x}=\omega_{F/F_{0}}$ et $\Delta_{\theta_{1}\otimes \theta_{2}} (\varpi)=\omega_{F/F_{0}}(\varpi(b_{1}-b_{2}))\lambda_{F_{\vert F_{0}}}(\psi_{0})^{-1}$.
\end{itemize}
Posons : $\Theta=\Delta_{\theta_{1}\otimes \theta_{2}}^{-1}\tl\theta_{1}\ol \Omega $ si $n>0$ et $\Theta =\mu^{-1}\Theta_{1}\ol \Theta_{2}$ si $n=0$.
Alors $\pi_{0}$ est la repr\'esentation de $GL_{2}(F_{0})$ associ\'ee au caract\`ere $\Theta$ par \cite{JL}.
\end{Corollaire}

\begin{proof} C'est une application des r\'esultats de \cite{BH3}, \S\S 19 et 34.
\end{proof}

\subsection{}  D'apr\`es la th\'eorie de Mackey, la restriction de $\pi_{0}$ \`a $GL_{2}(F_{0})^+$ est la somme des repr\'esentations $\pi_{0}^+$ et $\pi_{0}^-$. D'apr\`es \cite[p.738]{LL}, pour tout $x\in F^\x$ r\'egulier, 
\begin{equation*}
 \tr{\pi_{0}^+}(x)-\tr{\pi_{0}^-}(x)=\pm \lambda_{F_{\vert F_{0}}}(\psi_{0})\omega_{F/ F_{0}}\left(\frac{x-\ol x}{\varepsilon}\right) \frac{\Theta(x)+\Theta(\ol x)}{D(x)},
\end{equation*}
o\`u $D(x)= \left\vert \frac{ (x-\ol x)^2}{x\ol x}\right\vert_{F_{0}}^{1/2}$.
Mais, si $h=(h_{1}, h_{2})\in H$ est $G$-r\'egulier et $a\in F^\x$ tel que $a\ol a=h_{1}h_{2}$, $\ol a h_{1}$ est un \'el\'ement de $F^\x$ r\'egulier dont le conjugu\'e est $\ol a h_{2}$ et l'on a : 
\begin{equation*}
\begin{aligned}
&\tr{\pi}(h)=\tr{\pi}(\Phi^{-1}(i(\ol a h_{1})\cdot \ol a^{-1}))= \tr{\pi_{0}^+}(i(\ol a h_{1}))\Omega(\ol a^{-1}),\\
\text{et \quad } &\tr{{}^{g_{0}}\pi}(h)= \tr{\pi_{0}^-}(i(\ol a h_{1}))\Omega(\ol a^{-1}).
 \end{aligned}
\end{equation*}
Par suite (on note $\Delta=\Delta_{\theta_{1}\otimes \theta_{2}}$) :
\begin{equation*}
\begin{aligned}
\tr{\pi}(h)-\tr{{}^{g_{0}}\pi}(h)&= \pm  \lambda_{F_{\vert F_{0}}}(\psi_{0})\omega_{F/ F_{0}}\left(\frac{\ol ah_{1}-\ol ah_{2}}{\varepsilon}\right) \Omega(\ol a ^{-1}) \frac{\Theta(\ol ah_{1})+\Theta(\ol ah_{2})}{\Delta(\ol ah_{1})}\\
&=  \pm  \lambda_{F_{\vert F_{0}}}(\psi_{0})\mu(h_{1}-h_{2})\mu\left(\frac{\ol a}{\varepsilon}\right)\Delta^{-1}(\ol a)\Omega(\ol a ^{-1}) \\
&\quad \quad  \quad \quad \cdot
\frac{\Delta^{-1}(h_{1})\tl \theta_{1}(\ol ah_{1})\Omega(\ol ah_{2})+\Delta^{-1}(h_{2})\tl \theta_{1}(\ol ah_{2})\Omega(\ol ah_{1})}{D_{G}(h)}\\
&=\pm \lambda_{F_{\vert F_{0}}}(\psi_{0})\mu(\varepsilon)^{-1}\frac{\iota(h)}{D_{G}(h)}\mu(\ol a)\Delta^{-1} (\ol a)\\
&\quad \quad  \quad \quad \quad \quad  \quad  \quad      
 \cdot \left(\mu^{-1}\theta_{1}\otimes \theta_{2}(h)+\mu^{-1}\theta_{1}\otimes \theta_{2}(whw^{-1})\right)\\
 \end{aligned}
\end{equation*}
en remarquant que : $\Delta (x)=\mu(x)$ pour tout $x\in F_{\vert F_{0}}^1$ et $\tl \theta_{1}(\ol ah_{1})=\theta_{1}(h_{1}h_{2}^{-1})$ tandis que $\Omega(\ol a^{-1})\Omega (\ol ah_{2})=\Omega(h_{1}h_{2})\Omega(h_{1}^{-1})=\theta_{1}(h_{2})\theta_{2}(h_{2})$ car $h_{2}\in F_{\vert F_{0}}^1$.
De plus, $\lambda_{F_{\vert F_{0}}}(\psi_{0})\mu(\varepsilon)^{-1}$ est \'egal \`a $-1$ si $F$ est non ramifi\'ee et \`a 1 si $F$ est ramifi\'ee (par choix de $\varepsilon$ et $\mu$) et le caract\`ere de $F^\x$, $ a\mapsto \mu(\ol a)\Delta^{-1} (\ol a)$, est trivial si $F$ n'est pas ramifi\'ee et \'egal \`a $(\omega_{F/F_{0}}(\varpi(b_{1}-b_{2}))^{\val(a)}$ si $F$ est ramifi\'ee. Dans ce dernier cas, $a$ est de valuation paire si et seulement si $h_{1}h_{2}\in (1+{\goth p})_{F_{0}}^1$. 
En comparant \`a (\ref{IC-tr}), on conclut :

\begin{Proposition} Soit $\Pi=\{\pi, {}^{g_{0}}\pi\}$ un paquet endoscopique cuspidal de $G$ d\'ecrit en \ref{donnees} (dont on reprend les notations). On d\'efinit un  caract\`ere $\delta _{\theta_{1}\otimes \theta_{2}}$  de $F_{\vert F_{0}}^1$ par : $\delta _{\theta_{1}\otimes \theta_{2}}$ est trivial si $F$ n'est pas ramifi\'ee sur $F_{0}$ ; $\delta _{\theta_{1}\otimes \theta_{2}}$ est trivial sur $(1+{\goth p})^1_{\vert F_{0}}$ et $\delta _{\theta_{1}\otimes \theta_{2}}(-1)= \omega_{F/F_{0}}(\varpi (b_{1}-b_{2}))$ si $F$ est ramifi\'ee.\\
Alors le paquet $\Pi$ est l'image par l'application de transfert  du caract\`ere $\theta$ de $H$ d\'efini par :
\begin{equation*}
\forall h\in H, \quad \theta (h)=
\delta _{ _{\theta_{1}\otimes \theta_{2}}}(\Det h)\cdot \mu^{-1}\theta_{1}\otimes\theta_{2}(h). \end{equation*}
\end{Proposition}

\medskip	

En application de la proposition 11.4.1(a) de \cite{Ro}, on obtient :

\begin{Corollaire}\label{BC-endoscopiques} L'image par le changement de base stable du paquet $\Pi$ d\'ecrit en \ref{donnees} (dont on reprend les notations) est $\tl \pi =\indu_{\tl P}^{GL_{2}(F)} \mu(\theta)$ o\`u $\tl P$ est un sous-groupe parabolique de $GL_{2}(F)$ de facteur de Levi $H(F)$ et $\mu(\theta)$ le caract\`ere de $H(F)$ d\'efini par : 
\begin{equation*}
\forall h\in H(F), \quad\mu(\theta) (h)=
\mu^{-1}\tl \delta _{ _{\theta_{1}\otimes \theta_{2}}}(\Det h)\cdot \tl \mu^{-1}\tl \theta_{1}\otimes\tl \theta_{2}(h). \end{equation*}
Lorsque $F$ est ramifi\'ee sur $F_{0}$ et que $\Pi$ est l'unique paquet endoscopique de niveau $0$, l'image de $\Pi$ est $\tl \pi =\indu_{\tl P}^{GL_{2}(F)} \mu^{-1}\otimes \tl \chi \mu^{-1}$.
\end{Corollaire}

\section{Changement de base stable des paquets cuspidaux singletons.}\label{paquetssing} On proc\`ede en trois \'etapes :  la premi\`ere est consacr\'ee \`a la construction de repr\'esentations irr\'eductibles cuspidales de $\tl G=GL_{2}(F)$, $\tau$-invariantes et de caract\`ere central trivial sur $F_{0}^\x$ \`a partir de repr\'esentations tr\`es cuspidales de $G$. Par \cite[prop. 11.4.1(c)]{Ro}, ces repr\'esentations appartiennent \`a l'image d'un des deux changements de base, le ``stable'' ou le ``labile''. 
Dans la deuxi\`eme \'etape, on distingue parmi les repr\'esentations construites celles qui appar\-tien\-nent \`a l'image du changement de base stable. La derni\`ere \'etape d\'ecrit les changements de base stable et labile des paquets singletons de $G$.

\subsection{}\label{donnees2} Soit $(J,\lambda)$ un type simple maximal de $G$, c'est-\`a-dire :
\begin{enumerate}
\item[(a)]  {\it Cas de niveau 0 ($F$ ramifi\'ee sur $F_{0}$}) : $J$ est le sous-groupe parahorique maximal de $G$, c'est-\`a-dire $J=U_{0}(\scr L)$ avec $\scr L$ la cha\^\i ne autoduale de $F^2$, de p\'eriode 1 et d'invariant impair. La repr\'esentation  $\lambda$ est une repr\'esentation irr\'eductible de $J$, triviale sur le sous-groupe pro-unipotent $J_{1}$ de $J$ et dont la factorisation $\ol \lambda$ par $J/J_{1}\simeq SL_{2}(k_{0})$ est une repr\'esentation cuspidale de dimension $q-1$. Elle est associ\'ee \`a un caract\`ere $\theta$ r\'egulier d'ordre diff\'erent de 2 du groupe des \'el\'ements de norme 1 de l'extension quadratique $\ell$ de $k_{0}$ par :  
\item[]  $\quad \quad  \left\{ \begin{array}{ll}
\tr \ol \lambda(x)= (q-1)\theta(x) & \text{ si } x\in \{\pm 1\}\\ 
\tr \ol \lambda(xn)=-\theta(x) & \text{ si } x\in \{\pm 1\}, n\in N(k_{0})-\{\id\}\\ 
\tr \ol \lambda(x)=-(\theta(x)+\theta(\gamma(x)))& \text{ si } x\in\ell_{\vert k_{0}}^1 - \{\pm 1\}
\end{array}\right.$
\item[] o\`u $\gamma$ est l'\'el\'ement non trivial de $\Gal (\ell/k_{0})$ \cite{Sp}. 
\item[(b)] {\it Cas de niveau strictement positif} : $(J,\lambda)$ provient d'une strate gauche tr\`es cuspidale $(\scr L, n,n-1,b)$ avec $n>0$, c'est-\`a-dire, en notant $E=F[b]$ l'extension de $F$ engendr\'ee par $b$ et $L=E^\sigma$ la sous-extension de $E$ form\'ee des points fixes par $\sigma$,
\begin{equation*}\begin{split}
J={\goth o}^1_{E|L}U_{[{n+1\over 2}]}({\scr L})\supset &
J_1=(1+{\goth p}_E)^1_{|L}U_{[{n+1\over 2}]}({\scr L})\\
&\quad \quad \quad \quad \quad 
\supset H_1=(1+{\goth p}_E)^1_{|L}U_{[{n\over 2}]+1}({\scr L})
\end{split}
\end{equation*}
et $\lambda$ est obtenue \`a partir d'un caract\`ere $\theta$ de $H_1$ prolongeant le carac\-t\`ere $\psi_{b}$ de $U_{[\frac{n}{2}]+1}(\scr L)$, comme un  prolongement de l'unique repr\'esen\-ta\-tion irr\'eductible $\eta_{\theta}$ de $J_1$ contenant $\theta$ \cite[annexe]{Bl1}. On note $\omega_{\lambda}$ le caract\`ere central de $\lambda$.
\end{enumerate}
On note $(\pi, \scr V)$ l'induite compacte de $J$ \`a $G$ de $\lambda$. C'est une repr\'esentation irr\'eductible tr\`es cuspidale de $G$ (et toute repr\'esentation irr\'eductible tr\`es cuspidale de $G$ s'obtient ainsi).

\subsection{Construction dans le cas de niveau 0.}\label{construction0} On consid\`ere le type simple $(J,\lambda)$ d\'ecrit en \ref{donnees2} (a) (dont on reprend les notations) et on construit deux types simples maximaux $\tau$-invariants de $\tl G$, $(\tl J, \tl \Lambda)$ et $(\tl J, \tl \Lambda')$.

On choisit $\tl J= F^\x \tl U_{0}(\scr L)=\varpi^{\Bbb Z} \tl U_{0}(\scr L)$. La repr\'esentation $\tl \Lambda$ est un prolongement d'une repr\'esentation $\tau$-invariante $\tl \lambda$ de $\tl U_{0}(\scr L)$ triviale sur $\tl U_{1}(\scr L)$ qui se factorise en une repr\'esentation cuspidale $\ol {\tl \lambda}$ de  $U_{0}(\scr L)/U_{1}(\scr L)\simeq GL_{2}(k_{0})$. On d\'efinit donc ${\tl \lambda}$ comme suit.\\
On consid\`ere le caract\`ere  $\tl \theta$ de $\ell^\x$ d\'efini par : $\tl \theta(x)=\theta(x\gamma(x)^{-1})$, $x\in \ell^\x$. Il est  r\'egulier donc associ\'e \`a une repr\'esentation cuspidale $\ol {\tl \lambda}(\tl \theta)$ de $GL_{2}(k_{0})$, caract\'eris\'ee par  \cite[(6.4.1)]{BH3} :
\begin{equation*}
\left\{ \begin{array}{ll}
\tr \ol {\tl \lambda}(\tl \theta)(x)= (q-1)\tl \theta(x) & \text{ si } x\in k_{0}^\x\\ 
\tr \ol {\tl \lambda}(\tl \theta)(xn)=-\tl \theta(x) & \text{ si } x\in k_{0}^\x, n\in N(k_{0})-\{\id\}\\ 
\tr \ol {\tl \lambda}(\tl \theta)(x)=-(\tl \theta(x)+\tl\theta(\gamma(x)))& \text{ si } x\in\ell^\x - k_{0}^\x.
\end{array}\right.
\end{equation*}
 Alors $\tl \lambda$ est le rel\`evement de $\ol{\tl \lambda}(\tl \theta)$ \`a $\tl U_{0}(\scr L)$. Elle est de  caract\`ere central trivial sur ${\goth o}_{0}^\x$ et est $\tau$-invariante : 
\begin{equation}\label{tau0}
\forall g\in \tl U_{0}(\scr L),\quad  \tl \lambda (\tau(g))=\tl \lambda((\Det g)^{-1}g)=\tl \lambda(g)
\end{equation}
Elle poss\`ede deux prolongements \`a $\tl J$ $\tau$-invariants, $\tl \Lambda$ et $\tl \Lambda'$, qui diff\`erent par leur valeur en $\varpi$ : $\tl \Lambda (\varpi)= \omega_{\lambda}(-1)=-\tl \Lambda'(\varpi)$. 

On note alors
\begin{equation}\label{construction-choix0}
\tl \pi =\ind_{\tl J}^{\tl G}\tl \Lambda \quad \text{ et } \quad \tl \pi' =\ind_{\tl J}^{\tl G}\tl \Lambda'
\end{equation}
On prolonge $\tl \pi$ et $\tl \pi'$ \`a $\tl G\Gamma$ en imposant $\tl \Lambda(\tau)=\tl \Lambda'(\tau)=\tl \lambda(\tau)=1$ (\ref{tau0}). En notant $E$ l'extension quadratique non ramifi\'ee de $F$, dont le groupe multiplicatif $E^\x$ est plong\'e dans le normalisateur de ${\goth a}_{0}(\scr L)$, on a imm\'ediatement~:
\begin{equation}\label{tr-lambda-0}
\forall x\in E^\x, x\not \in F^\x, \quad \tr \tl \Lambda(x)=\tr \lambda(\N(x)).
\end{equation}

\subsection{Construction dans le cas de niveau strictement positif.}\label{construction}
C'est l'analogue de \cite{Bl2} dans le cas de dimension deux. 

On consid\`ere les sous-groupes ouverts compacts modulo le centre de $\widetilde G$, $\Gamma$-invariants :
$$\aligned
&\widetilde H_1=(1+{\goth p}_E)\widetilde U_{[{n\over 2}]+1}({\scr L}),
\quad 
\widetilde J_1=(1+{\goth p}_E)\widetilde U_{[{n+1\over 2}]}({\scr L}),
\cr
&\widetilde J_c={\goth o}_E^\times \widetilde U_{[{n+1\over 2}]}({\scr
L}), \quad
\widetilde J=E^\times \widetilde U_{[{n+1\over 2}]}({\scr L})=E^\times
\widetilde J_c.\endaligned$$
Sur $\widetilde H_1$, on consid\`ere le caract\`ere $\tl \theta=\theta\circ \Nc$ \cite[cor. 3.2]{Bl2}. Si $n$ est pair, il existe une unique repr\'esentation $\tl \eta$ de $\tl J_{1}$ contenant $\tl \theta$, n\'ecessairement $\tau$-invariante. Si $n$ est impair, on pose $\tl \eta=\tl \theta$.\\
Puisque la dimension de $\tl \eta$ et l'ordre de $\tl J_{c}/\tl J_{1}$ sont premiers entre eux, il exis\-te des prolongements de $\tl \eta$, et en particulier des prolongements $\tau$-invariants de caract\`ere central trivial sur ${\goth o}_{0}^\x$. Deux tels prolongements ont m\^eme res\-triction au sous-groupe $\tl J_{0}={\goth o}_{0}^\x N_{E_{\vert L}}({\goth o}_{E}^\x)\tl J_{1}$. On distingue alors trois cas~:
\begin{itemize}
\item [(nr-nr)] {\it $F$ n'est pas ramifi\'ee sur $F_{0}$} : l'extension $E$ est alors ramifi\'ee et $\tl \eta$ est de dimension 1. De plus, $\mu$ est d'ordre 2 et les deux changements de base de $\pi$ ont m\^eme caract\`ere central \'egal \`a $\tl{\omega_{\lambda}}$.  On impose donc que le caract\`ere central de $\tl \Lambda$ soit $\tl{\omega_{\lambda}}$. Il existe alors deux prolongements de $ \tl{\omega_{\lambda}}\tl \eta$ \`a $\tl J$~: ils ont m\^eme restriction $\tl \lambda$ \`a $\tl J_{c}$ et sont \'egaux \`a $\pm 1$ en $\varpi_{L}$.
\item[(r-r)] {\it  $F$ est ramifi\'ee sur $F_{0}$ et $n$ est impair} : l'extension $E$ n'est pas ramifi\'ee sur $F$ et $\varpi^{-n}b$ est un \'el\'ement de ${\goth o}_{E}^\x$ invariant par $\sigma$ : $L=F_{0}[\varpi^{-n}b]$ n'est pas ramifi\'ee sur $F_{0}$ et $E$ est ramifi\'ee sur $L$. Le groupe $\tl J_{0}$ est d'indice 2 dans $\tl J_{c}$ et $\tl J=F^\x\tl J_{c}$. Il existe quatre prolongements $\tau$-invariants de $\tl \eta=\tl \theta$ \`a $\tl J$ de caract\`ere central trivial sur $F_{0}^\x$, deux de caract\`ere central $\tl{\omega_{\lambda}}$ et deux de caract\`ere central $(-1)^{\val}\tl{\omega_{\lambda}}$. Mais, si l'on consid\`ere l'autre prolongement de $\eta$ \`a $J$, on retrouve ces m\^emes prolongements de $\tl \eta$. 
\item[(r-nr)] {\it $F$ est ramifi\'ee sur $F_{0}$ et $n$ est pair} : l'extension $E$ n'est pas ramifi\'ee sur $F$ et $L$ est engendr\'ee sur $F_{0}$ par $\varpi^{-n+1}b$, une uniformisante de $E$. Donc $E$ n'est pas ramifi\'ee sur $L$. Le groupe $\tl J_{0}$ est d'indice $q+1$ dans $\tl J_{c}$ et $\tl J=F^\x\tl J_{c}$. Il existe $2(q+1)$ prolongements $\tau$-invariants de $\tl \eta$ \`a $\tl J$  de caract\`ere central trivial sur $F_{0}^\x$, une moiti\'e  de caract\`ere central $\tl{\omega_{\lambda}}$ et l'autre de caract\`ere central $(-1)^{\val}\tl{\omega_{\lambda}}$. Mais, si l'on consid\`ere les $(q+1)$ prolongements de $\eta$ \`a $J$ on retrouve les m\^emes prolongements de $\tl \eta$. 
\end{itemize}

\begin{Remarque} Ces trois cas, auxquels on ne cesse de faire r\'ef\'erence par la suite, se rep\`erent par la ramification de $F$ sur $F_{0}$ puis celle de $E$ sur $E^\sigma$. Les deux premiers cas correspondent \`a des paires $(E/F,\sigma)$ ``paire'' au sens de la d\'efinition 1.2 de \cite{HM} tandis que le dernier cas correspond \`a une paire $(E/F,\sigma)$ ``impaire'', tout comme le cas de niveau 0 vu pr\'ec\'edemment. Cette diff\'erence se refl\`ete dans les \'enonc\'es des proposition \ref{stable} et th\'eor\`eme \ref{BCdim2}.
\end{Remarque}

\begin{Lemme}\label{construction-choix}
 Il existe une unique repr\'esentation $\tl \Lambda$ de $\tl J\cdot \Gamma$ v\'erifiant :
\begin{itemize}
\item[(i)] $\tl \Lambda_{\vert \tl J}$ est un prolongement $\tau$-invariant de  $\tl \eta$, de caract\`ere central $\tl{\omega_{\lambda}}$, 
\item[(ii)] pour tout $x\in E^\x$, $\tr \tl \Lambda (x\tau)=\tr \lambda(\Nc (x))$. 
\end{itemize}
\end{Lemme}

On note alors
\begin{equation}\label{construction-choix2}
\tl \pi =\ind_{\tl J}^{\tl G}\tl \Lambda_{\vert \tl J}
\end{equation}
que l'on prolonge \`a $\tl G\Gamma$ par $\ind_{\tl J \Gamma}^{\tl G \Gamma}\tl \Lambda$.

\begin{proof}
On rappelle \cite[lemme 3.3]{Bl2} qu'il existe un unique prolongement de $\tl \eta$ \`a $\tl J_{1}\Gamma$ tel que : 
\begin{equation*}
\forall g\in \tl J_1, \tr \tl \eta (g\tau)= \tr \eta (\Nc (g)).
\end{equation*}
Soit $\tl \Lambda$ un prolongement de $\tl \eta$ de caract\`ere central $\tl{\omega_{\lambda}}$. On le prolonge \`a $\tl J\cdot \Gamma$ en imposant : $\tl \Lambda(\tau)=\tl \eta (\tau)$. Dans les cas (nr-nr) et (r-r), $\tl \Lambda$ est de dimension 1 et la condition (ii) du lemme \'equivaut \`a : $\tl \Lambda(\varpi_{E})=1$ dans le cas (nr-nr) ; $ \tl \Lambda(\zeta)=1$ o\`u $\zeta$ est une racine primitive $(q^2-1)$-i\`eme de $1$ contenue dans $E$  dans le cas (r-r). Ceci d\'efinit un et un seul prolongement de la liste.\\
Dans le cas (r-nr), seule la restriction $\tl \lambda$ de $\tl \Lambda$ \`a $\tl J_{c}$ importe. Le calcul de la trace de $\lambda$ sur les \'el\'ements de ${\goth o}_{E_{\vert L}}^1$ est effectu\'e au paragraphe \ref{Traces} (\ref{precisions}, (2)). Il existe donc un unique caract\`ere $\xi$ de ${\goth o}_{E_{\vert L}}^1H_{1}$, prolongeant $\theta$ tel que :
\begin{equation}\label{tr-lambda}
\forall x\in {\goth o}_{E_{\vert L}}^1, x\not \in{\goth o}_{\vert F_{0}}^1, \tr \lambda(x)=\epsilon \xi(x),
\end{equation}
o\`u $\epsilon$ est \'egal \`a -1 si $\dim \eta =q$, 1 sinon.\\
Pour le calcul de $\tr \tl \lambda$, on distingue deux cas suivant la dimension de $\tl \lambda$. Si $\tl \lambda$ est de dimension 1, $\lambda$ l'est \'egalement et le r\'esultat est imm\'ediat. Supposons donc $\tl \lambda$ de dimension $q$. On introduit  alors le groupe ${\scr X}$ des caract\`eres de $\tl J_{c}/{\goth o}^\x\tl J_{1}$ qui  s'identifie au groupe des caract\`eres de $k_{E}^\x$ triviaux sur $k^\x$, groupe cyclique d'ordre $q+1$ engendr\'e par un caract\`ere $\tau$-invariant $\tl \kappa$. Gr\^ace aux calculs effectu\'es en \ref{Tr32} avec les notations pr\'ecis\'ees au (3) de \ref{precisions}, on est assur\'e de l'existence de $q+1$ entiers $m_{i}$ de somme $q$ tels que :
$$\forall x\in {\goth o}_{E}^\x, \tr \tl \lambda(x)=\oplus_{i=0}^q m_{i}(\tl \kappa^i\tl \xi)(x).$$
Chaque composante isotypique \'etant $\tau$-invariante, il existe $q+1$ entiers $n_{i}$ tels que :
$$\forall x\in {\goth o}_{E}^\x, \tr \tl \lambda(x\tau)=\oplus_{i=0}^q n_{i}(\tl \kappa^i\tl \xi)(x)\quad \text{ et } \quad \sum_{i=0}^q n_{i}=\tr \tl \eta (\tau)=\dim \eta.$$ 
On raisonne alors comme dans la d\'emonstration de la proposition 3.5 de \cite{Bl2}, et on obtient :
\begin{itemize}
\item[-] si $\dim \eta =1$, tous les entiers $n_{i}$ sont nuls sauf un \'egal \`a $1$ ;
\item[-] si $\dim \eta =q$, tous les entiers $n_{i}$ sont \'egaux \`a $1$ sauf un qui est nul. 
\end{itemize}
Il existe donc un et un seul prolongement $\tl \lambda$ de $\tl \eta$ \`a $\tl J_{c}$ tel que :
$$\forall x\in {\goth o}_{E}^\x, x\not \in{\goth o}^\x, \tr \tl \lambda(x\tau)=\varepsilon \tl \xi(x)=\tr \lambda(\Nc (x)).$$
Ceci termine la d\'emonstration du lemme.
\end{proof}
 
\subsection{Stabilit\'e.}
\begin{Proposition}\label{stable} Soit $\tl \pi=\ind_{\tl J}^{\tl G} \tl \Lambda$ la repr\'esentation de $GL_{2}(F)$ d\'efinie en $(\ref{construction-choix0})$ ou $( \ref{construction-choix2})$. \\
(i) Elle appartient \`a l'image du changement de base stable si et seulement si 
\begin{enumerate}
\item[-] elle est de niveau $0$ ou
\item[-] elle est de niveau strictement positif et l'extension $E$ associ\'ee est non ramifi\'ee sur $F$ et sur $L:=E^\sigma$.
\end{enumerate}
(ii) Lorsque $\tl \pi$ appartient \`a l'image du changement de base labile, $\tl \pi \cdot \mu^{-1}\circ \Det$ appartient \`a l'image du changement de base stable.
\end{Proposition}

\begin{proof} L'assertion (ii) n'est qu'une redite de \cite[\S 11.4]{Ro}. Pour l'as\-ser\-tion (i), on distingue deux cas.\\
On suppose d'abord que $\tl \Lambda$ est de dimension 1. Alors par \cite[\S 19]{BH3} $\tl \pi$ est associ\'ee \`a la paire admissible $(E, \tl \Lambda_{\vert E^\x})$. Mais $E$ est soit d'indice de ramification 2 sur $F$ , soit non ramifi\'ee sur $F$ mais ramifi\'ee sur $L$. Autrement dit la paire $(E/F,\sigma)$ est paire au sens de \cite[d\'ef. 1.2]{HM}. Comme $\tl \Lambda$ est triviale sur $L^\x$ (lemme \ref{construction-choix} (ii)), $\tl \pi$ est $GL_{2}(F_{0})$-distingu\'ee \cite[Th\'eor\`eme 1.1]{HM} donc $\tl \pi$ appartient au changement de base labile \cite[Th\'eor\`eme 7]{Fl2}.\\
On suppose maintenant que $\tl \Lambda$ est de dimension au moins 2. Il n'y a plus de rapport simple entre la trace et la trace tordue de $\tl \Lambda$ ce qui rend l'emploi du crit\`ere pr\'ec\'edent moins adapt\'e que celui de la d\'efinition. Gr\^ace \`a \cite{Ro}, il n'y a que deux possibilit\'es : soit le caract\`ere tordu de $\tl \pi$, soit celui de $\tl \pi\cdot \mu\circ \Det^{-1}$ est constant sur les classes de $\tau$-conjugaison stable de $\tl G$. On peut donc d\'eduire si $\tl \pi$ appartient \`a l'image du changement de base stable en \'evaluant son caract\`ere tordu en deux \'el\'ements $g$ et $g'$ stablement $\tau$-conjugu\'es et  tels que $\tr \tl \pi(g\tau)\neq 0$ et $\mu\circ \Det (g^{-1}g')\neq 1$. \\
On choisit $g=\zeta$ ($\zeta$ est une racine primitive $(q^2-1)$-i\`eme de 1 contenue dans ${\goth o}_{E}^\x$) et $g'=\varpi_{L}g$ o\`u $\varpi_{L}=\zeta^{\frac{q+1}{2}}\varpi$. 
Tous deux sont de norme cyclique $x=\zeta^{1-q}$ et $\mu\circ \Det (\varpi_{L})=\mu\circ N_{L_{\vert F_{0}}}(\varpi_{L})=-1$. Notons que $x\in \tl U_{0}(\scr L)$ est tr\`es r\'egulier :  \quad $\forall h\in \tl G, \quad h^{-1}xh\in \tl U_{0}(\scr L) \impl h\in F^\x \tl U_{0}(\scr L)$.\\
Par cons\'equent, si $\tl \pi$ est de niveau 0, la formule de Mackey donne :
\begin{equation*}\begin{split}
&\tr \tl \pi(\zeta\tau)=\sum_{\begin{subarray}{c} h\in \tl G/\tl J\\ h^{-1}\zeta\tau(h)\in \tl J\end{subarray}} \tr \tl \Lambda (h^{-1}\zeta \tau(h)\tau)= \tr \tl \Lambda (\zeta\tau)=\tr \lambda(x) \\
\text{ et de m\^eme, } & \tr \tl \pi(\zeta\varpi_{L}\tau)=\omega_{\lambda}(-1)\tr \tl \pi(\zeta^\frac{q+3}{2}\tau)=\omega_{\lambda}(-1)\tr \lambda(-x)=\tr \tl \pi(\zeta\tau).
\end{split}\end{equation*}
D'apr\`es la d\'efinition de $\lambda$ (cf. \ref{donnees2} (a)), $\tr \lambda(x)$ n'est pas nul d\`es que le caract\`ere $\theta$ associ\'e n'est pas d'ordre 4. Si $\theta$ est d'ordre 4 ($q\equiv -1 \mod 4$), on refait le m\^eme raisonnement avec $g=\zeta^2$  sachant que $\tr \tl \pi(\zeta^2\tau)\neq0$. \\
Quand $\tl \pi$ est de niveau 0, elle appartient \`a l'image du changement de base stable.

Supposons maintenant que $\tl \pi$ est de niveau strictement positif. Notons $\tl K(\scr L)$ le normalisateur dans $\tl G$ de l'ordre associ\'e \`a $\scr L$, $\tl K(\scr L)=F^\x\tl U_{0}(\scr L)$, et $\tl \sigma=\ind_{\tl J\Gamma}^{\tl K(\scr L)\Gamma}\tl \Lambda$. Alors $\tl \pi$ est $\ind_{\tl K(\scr L)\Gamma}^{\tl G\Gamma}\tl \sigma$ et en appliquant deux fois la formule de Mackey, on a :
\begin{equation}\label{tr-tlpi(iii)}
\begin{split}
&\tr \tl \pi(\zeta\tau)=\tr \tl \sigma(\zeta\tau)=\sum_{\begin{subarray}{c} h\in \tl U_{0}(\scr L)/\tl J_{c}\\ h^{-1}\zeta\tau(h)\in \tl J_{c}\end{subarray}} \tr \tl \Lambda (h^{-1}\zeta \tau(h)\tau)\\ \text{ et } \quad 
&\tr \tl \pi(\zeta\varpi_{L}\tau)=\omega_{\lambda}(-1)\sum_{\begin{subarray}{c} h\in \tl U_{0}(\scr L)/\tl J_{c}\\ h^{-1}\zeta^\frac{q+3}{2}\tau(h)\in \tl J_{c}\end{subarray}} \tr \tl \Lambda (h^{-1}\zeta^\frac{q+3}{2} \tau(h)\tau)
\end{split}\end{equation}
On reprend alors le raisonnement de \cite[\S 4.7]{Bl2} en remarquant que $\Gal (E/F_{0})$ est ab\'elien (car d'ordre 4). On obtient :
\begin{equation*}\begin{split}
\tr \tl \pi(\zeta\tau)&=\sum_{\gamma\in \Gal (E/F)}\tr\tl \lambda(\gamma(\zeta)\tau) =\tr \lambda(x)+\tr \lambda(x^{-1})\\&=\omega_{\lambda}(-1)(\tr \lambda(-x)+\tr\lambda(-x^{-1}))=\tr \tl \pi(\zeta\varpi_{L}\tau)
\end{split}\end{equation*}
Or $\lambda$ est associ\'e \`a un caract\`ere $\xi$ par (\ref{tr-lambda}) : ou bien $\xi(x)+\xi(x^{-1})\neq 0$, ou bien  $\xi(x^2)+\xi(x^{-2})\neq 0$. On proc\`ede comme dans le cas de niveau 0. 
\end{proof}

\subsection{Description des changements de base.}
\begin{Theoreme}\label{BCdim2} Soit $(\pi, {\scr V})$ une repr\'esentation irr\'eductible tr\`es cuspidale de $G$ et $(J,\lambda)$ un type simple maximal d\'efinissant $\pi$ comme au paragraphe $\ref{donnees2}$, (a) ou (b). On note $\tl \pi=\ind_{\tl J}^{\tl G}\tl \Lambda$ la repr\'esentation irr\'eductible cuspidale de $\tl G$ d\'efinie en $(\ref{construction-choix0})$ ou $(\ref{construction-choix2})$ selon que $\pi$ est de niveau $0$ ou non. \\
(i) {\rm $F$ n'est pas ramifi\'ee sur $F_{0}$ :} l'image de $\pi$ par le changement de base stable est $\tl \pi\cdot \mu^{-1}\!\circ\! \Det $ tandis que celle par le changement de base labile est $\tl \pi$.\\
(ii) {\rm $F$ est ramifi\'ee sur $F_{0}$ et $E$ est ramifi\'ee sur $E^\sigma$ :} on note $\chi$ le caract\`ere d'ordre $2$ de $J$ trivial sur $J_{1}$ et on pose : $\pi'= \ind_{J}^G \lambda\cdot \chi\;$ et $\; \tl \pi'=\ind_{\tl J}^{\tl G}\tl \Lambda\cdot \tl \chi$.\\
L'image de $\pi$ par le changement de base stable est $\tl \pi\cdot \mu^{-1}\circ \Det $ si $q\equiv 1 \mod 4$,\, $\tl \pi'\cdot \mu^{-1}\circ \Det $ sinon. L'image de $\tl \pi$ par le changement de base labile est alors $\tl \pi$ si $q\equiv 1 \mod 4$ et $\tl \pi'$ sinon.\\
(iii) {\rm $F$ est ramifi\'ee sur $F_{0}$ et  $E$ n'est pas ramifi\'ee sur $E^\sigma$ (y compris cas de niveau $0$) :} l'image de $\pi$ par le changement de base stable est $\tl \pi$, celle par le changement de base labile $\tl \pi\cdot \mu\circ \Det$.
\end{Theoreme}

Suit la d\'emonstration de ce th\'eor\`eme.

\subsection{} Gr\^ace \`a \cite[\S 11.4]{Ro} il suffit d'\'etablir le r\'esultat pour le changement de base stable.
Notons $\tl \pi_{st}$ la repr\'esentation de $\tl G$ pr\'essentie \^etre l'image de $\pi$ par le changement de base stable dans le th\'eor\`eme et montrons qu'elle v\'erifie bien l'identit\'e de caract\`eres d\'ecrivant ce dernier \cite[\S\S\,  4.11, 12.5 et 11.4]{Ro}, c'est-\`a-dire : 
pour tout $g\in \tl G$ dont la
norme cyclique $\Nc (g)$ est r\'eguli\`ere elliptique et tout $x\in \Nc (g)$,
\begin{equation}\label{id}
\tr \tl \pi_{st}  (g\tau)=c_\tau(\tl \pi_{st} ) \tr
\pi (x),
\end{equation}
o\`u $\tl \pi_{st} (\tau)$ est un op\'erateur d'entrelacement d'ordre 2
entre $\tl \pi_{st} $ et $\tl \pi_{st} ^\tau$ et $c_\tau(\tl \pi_{st} )$ un signe ne
d\'ependant que du choix de $\tl \pi_{st} (\tau)$. Dans la construction pr\'ec\'edente, on a choisi $\tl \pi_{st}(\tau)$ pour que $c_\tau(\tl \pi_{st} )$ soit \'egal \`a 1.

\smallskip	

Soit $g\in \tl G$. On suppose que $x=\N(g)$ est un \'el\'ement de $G$, r\'egulier et elliptique. On note $T(F_{0})$ le centralisateur de $x$ dans $G$. C'est un tore compact de $G$ et $T(F)$ est $\tau$-invariant et isomorphe soit au groupe multiplicatif d'une extension quadratique $E_{x}$ de $F$, soit \`a $F^\x \x F^\x$.\\
La formule de Mackey fournit une expression des traces de $\pi (x)$ et $\tl \pi (g\tau)$ :
\begin{equation}\label{tr1}
\tr \pi (x)= \sum_{\smallmatrix y\in G/J \cr y^{-1}xy\in
J\cr\endmatrix}
\tr \lambda(y^{-1}xy) \quad \text{ et } \quad
\tr \tl \pi_{st} (g\tau)=\sum_{\smallmatrix h\in \tl G/\tl J\cr
h^{-1}g\tau(h)\in \tl J\cr\endmatrix} \tr \tl \Lambda_{st} (h^{-1}g\tau(h)\tau).
\end{equation}

\subsection{} On suppose d'abord que $\pi$ est de niveau strictement positif. La d\'emons\-tra\-tion est semblable \`a celle du th\'eor\`eme 3.7 de \cite{Bl2}. \\
On introduit le sous-groupe distingu\'e de $\tl G$, not\'e $\tl G^+$, d\'efini comme le noyau du caract\`ere $\mu\circ \Det$. Pour tout sous-groupe $H$ de $\tl G$, on abr\`ege $H\cap \tl G^+$ en $H^+$.

\begin{Lemme} Soient $g\in \tl G$ et $x\in G$ comme ci-dessus. On note ${\scr Cl}_{\tau}^{st}(g)$ la  classe de $\tau$-conjugaison stable de $g$ et ${\scr Cl}^{st}(x)$ la classe de conjugaison stable de $x$.\\
(i) Le groupe $\tl G^+$ est r\'eunion de classes de $\tau$-conjugaison.\\
(ii) Si $\mu\circ \Det$ est d'ordre $2$, $\tl G^+$ rencontre la classe de $\tau$-conjugaison stable de $g$ et ${\scr Cl}_{\tau}^{st}(g)^+$ contient la moiti\'e des classes de $\tau$-conjugaison contenues dans $ {\scr Cl}_{\tau}^{st}(g)$.\\
(iii) Si $\mu\circ \Det$ est d'ordre $4$, ou bien ${\scr Cl}_{\tau}^{st}(g)^+$ n'est pas vide et  alors ${\scr Cl}_{\tau}^{st}(g)^+$ contient la moiti\'e des classes de $\tau$-conjugaison contenues dans $ {\scr Cl}_{\tau}^{st}(g)$ ; ou bien ${\scr Cl}_{\tau}^{st}(g)$ ne rencontre pas $\tl G^+$, et alors : 

\quad \quad \quad \quad \quad $ {\scr Cl}_{\tau}^{st}(g)\cap \tl J=\emptyset\quad $ et $\quad {\scr Cl}^{st}(x)\cap J=\emptyset$.\\
(iv) Chaque classe de $\tau$-conjugaison dans ${\scr Cl}_{\tau}^{st}(g)^+$ contient deux classes de $\tau$-$\tl G^+$-conjugaison.
\end{Lemme}

\begin{proof}
Pour la premi\`ere assertion, il suffit de remarquer que $\mu\circ \Det$ est constant sur les classes de $\tau$-conjugaison.\\
Pour la suite, on sait que les classes de $\tau$-conjugaison de $g$ dans $ {\scr Cl}_{\tau}^{st}(g)$ sont param\'etr\'ees par $H^1(\Gamma, T(F))$ : \`a un cocycle $c$, on associe la classe de $\tau$-conjugaison de $c(\tau)g$.\\
Si $T(F)\simeq E_{x}^\x$, il y a deux classes de $\tau$-conjugaison, celle de $g$ et celle de $ag$ o\`u $a\in L_{x}^\x, a\not \in N_{{E_{x}}_{\vert L_{x}}}(E_{x}^\x)$ et $L_{x}={E_{x}^\sigma}$.  Notons que $L_{x}$ est une extension quadratique de $F_{0}$, n\'ecessairement distincte de $F$. Mais alors, $\Det a$ appartient \`a $N_{{L_{x}}_{\vert F_{0}}}(L_{x}^\x)$ sans \^etre un carr\'e de $F_{0}$ donc $\mu\circ \Det (a)$ est \'egal \`a $-1$. Ainsi, $\mu\circ \Det$ prend deux valeurs oppos\'ees sur les deux classes de $\tau$-conjugaison.\\
Si $T(F)\simeq F^\x \x F^\x$, il y a quatre classes de $\tau$-conjugaison, celles de $g$, $(1, \alpha_{0})g$, $(\alpha_{0}, 1)g$ et $(\alpha_{0}, \alpha_{0})g$. Encore une fois $\mu\circ \Det$ prend deux valeurs oppos\'ees.\\
Par cons\'equent, si $\mu \circ \Det$ est d'ordre 2, la moiti\'e des classes de $\tau$-conjugaison dans $ {\scr Cl}_{\tau}^{st}(g)$ sont contenues dans $\tl G^+$.\\
Supposons que $\mu\circ \Det$ est d'ordre 4. Ou bien $\mu\circ\Det (g)\in \{ \pm 1\}$ et la moiti\'e des classes de $\tau$-conjugaison dans $ {\scr Cl}_{\tau}^{st}(g)$ sont contenues dans $\tl G^+$ ; ou bien $\mu\circ \Det (g)\in \{ \pm i\}$ et ${\scr Cl}_{\tau}^{st}(g)^+$ est vide. On termine la d\'emonstration en remarquant que : $\mu\circ \Det (\tl J)=\{ \pm 1\}$ et $\Det (J)\subset 1+{\goth p}$. En effet, l'hypoth\`ese sur $\Det g$ implique que celui-ci est de valuation impaire et par suite, que $\Det x$ appartient \`a $-1+{\goth p}$.\\
(iv) Si $g_{1}$ et $g_{2}$ sont deux \'el\'ements de $\tl G^+$ $\tau$-$\tl G$-conjugu\'es par deux \'el\'ements $h_{1}$ et $h_{2}$ alors $h_{1}h_{2}^{-1}$ appartient au $\tau$-centralisateur de $g_{1}$, en particulier la norme de son d\'eterminant est 1. Par cons\'equent, si $\mu$ est d'ordre 2, $h_{1}$ et $h_{2}$ diff\`erent d'un \'el\'ement de $\tl G^+$. Si $\mu$ est d'ordre 4, le $\tau$-centralisateur de $g_{1}$ contient un \'el\'ement $a$ tel que $\mu\circ \Det (a)=-1$ (cf. (ii) et (iii)) et $h_{1}$ et $h_{2}$ diff\`erent d'un \'el\'ement de $\tl G^+$ ou de $\tl G^+a$. Dans les deux cas, il y a deux classes de $\tau$-$\tl G^+$-conjugaison dans une classe de $\tau$-conjugaison.
\end{proof}

Dans le cas o\`u ${\scr Cl}_{\tau}^{st}(g)^+$ est vide, l'\'egalit\'e (\ref{id}) est satisfaite puisque les deux sommes de (\ref{tr1}) sont nulles. On peut donc supposer que ${\scr Cl}_{\tau}^{st}(g)^+$ n'est pas vide et, quitte \`a changer $g$ dans sa classe de $\tau$-conjugaison stable sans changer de $x$, que $g\in \tl G^+$. 

Notons ${\scr C}_{J}^{st}(x)$  l'ensemble des classes de $J$-conjugaison contenues dans ${\scr Cl}^{st}(x)\cap J$ et ${\scr C}_{\tau,\tl J^+}^{st}(g)$ celui des classes de $\tau$-$\tl J^+$-conjugaison contenues dans ${\scr Cl}_{\tau}^{st}(g)\cap \tl J^+$. Notons \'egalement $n(x)$ (resp. $\tl n(g)$) le nombre de classes de conjugaison (resp. $\tau$-conjugaison) contenues dans ${\scr Cl}^{st}(x)$ (resp. ${\scr Cl}_{\tau}^{st}(g)$).
En suivant le raisonnement de \cite[\S 4.2]{Bl2}, on a :
\begin{equation}\label{tr2}
\begin{aligned}
&\tr \pi (x)=\frac{1}{n(x)}\sum_{x'\in
{\scr C}_{J}^{st}(x)} c(x')\tr \lambda(x') \\ 
Ê\text{ et } \quad&
\tr \tl \pi_{st} (g\tau)=\frac{1}{\tl n(g)} \sum_{g'\in {\scr C}_{\tau,\tl J^+}^{st}(g)} \tl
c(g')\tr \tl \Lambda_{st}(g'\tau)
\end{aligned}
\end{equation}
o\`u $\quad c(x')=[T(F_0) : T(F_0)\cap yJy^{-1}] $  si $x'=y^{-1}xy,\,\, y\in G$\\
et $\quad\tl c(g')=[T(F_0) : T(F_0)\cap h\tl J{h}^{-1}] $ si $g'=h^{-1}g\tau (h),\,\, h\in \tl G$.

\subsection{} On introduit les sous-groupes $J'$ et $\tl J'$ de $J$ et $\tl J^+$ respectivement, en distinguant deux cas  :
\begin{itemize}
\item[cas 1 :] {\it $F/F_{0}$ est non ramifi\'ee ou  $q\equiv 1\mod 4$}. Dans ce cas, $J'=F_{\vert F_{0}}^1J_{1}$ et $\tl J'=F^\x\tl J_{1}$ ; 
\item[cas 2 :] {\it $F/F_{0}$ est ramifi\'ee et $q\equiv -1 \mod 4$}. Alors $J'=J_{1}$ et $\tl J'=F_{0}^\x\tl J_{1}$.
\end{itemize}
Dans les cas o\`u $J$ est plus grand que $J'$ ou $\tl J^+$ plus grand que $N_{E_{\vert L}}({\goth o}_{E}^\x)\tl J'$, c'est-\`a-dire lorsqu'on se trouve dans le cas (r-nr) ou le cas (r-r) et $q\equiv -1 \mod 4$, on consid\'ere les repr\'esentations par ``paquets'' fa\-bri\-qu\'es de la fa\c con suivante.

Du c\^ot\'e de $G$, on consid\`ere tous les prolongements \`a $J$ de la restriction \`a $J'$ de $\lambda$. Ils sont au nombre de $d=[J:J']$ et de la forme $\lambda \otimes \kappa^r$, $0\leq r\leq d-1$,  o\`u $\kappa$ est un caract\`ere de $J$, trivial sur $J'$ et d'ordre $d$ (qui s'identifie \`a un caract\`ere du groupe engendr\'e par $\zeta^{q-1}$ dans le cas (r-nr), -1 dans le cas (r-r)).\\
Notons  $\tl \kappa$ le caract\`ere de $\tl J$ que l'on obtient sur $\tl J$ en relevant le caract\`ere $\kappa\circ \N$ de $\tl J/N_{E_{\vert L}}({\goth o}_{E}^\x)\tl J'$. On remplace $\pi$ et $\tl \pi_{st}$ par :
\begin{equation*}
\begin{aligned}
\pi=\oplus_{r=0}^{d-1}\pi_{r} \quad& \text{ o\`u } \quad \pi_{r}=\ind_{J}^G \lambda\otimes \kappa^r \\
\tl \pi_{st}=\oplus_{r=0}^{d-1}\tl \pi_{r} \quad& \text{ o\`u } \quad \tl \pi_{r}=\ind_{\tl J}^{\tl G}\tl \Lambda_{st}\otimes\tl \kappa^r 
\end{aligned}
\end{equation*}
Alors, en sommant les expressions (\ref{tr2}) des caract\`eres des $\pi_{r}$ d'une part, et des $\tl \pi_{r}$ d'autre part, on obtient :
\begin{equation*}\begin{aligned}
&\tr \pi (x)=\frac{1}{n(x)}\sum_{x'\in
{\scr C}_{J}^{st}(x)\cap J'}d c(x')\tr \lambda(x')=\frac{1}{n(x)}\sum_{x'\in
{\scr C}_{J'}^{st}(x)}\frac{d}{d(x')} c(x')\tr \lambda(x')Ê\\
&
\begin{aligned}\text{ et} \quad\tr \tl \pi_{st} (g\tau)=
&\frac{1}{\tl n(g)} \sum_{g'\in {\scr C}_{\tau,\tl J^+}^{st}(g)\cap N_{E_{\vert L}}({\goth o}_{E}^\x)\tl J'}  d\tl c(g')\tr \tl \Lambda_{st}(g'\tau)\\ 
=&\frac{1}{\tl n(g)} \sum_{g'\in {\scr C}_{\tau,N_{E_{\vert L}}({\goth o}_{E}^\x)\tl J'}^{st}(g)}\frac{d}{\tl d(g')}\tl c(g')\tr \tl \Lambda_{st}(g'\tau)\end{aligned}\end{aligned}
\end{equation*}
o\`u $d(x')$ est le nombre de classes de $J'$-conjugai\-son dans l'intersection de la classe de $J$-conjugaison de $x'$ avec $J'$ et $\tl d(g')$ celui des classes  de $\tau$-$N_{E_{\vert L}}({\goth o}_{E}^\x)\tl J'$-conjugaison dans l'intersection de la classe de $\tau$-$\tl J^+$-conjugaison de $g'$ avec $N_{E_{\vert L}}({\goth o}_{E}^\x)\tl J'$. \\
Remarquons que les nombres $d(x')$ et $\tl d(g')$ sont invariants quand on conjugue $x'$ par un \'el\'ement de $J$, respectivement $\tau$-conjugue $g'$ par un \'el\'ement de $N_{E_{\vert L}}({\goth o}_{E}^\x)\tl J'$.

\smallskip	

\noindent On uniformise les notations en posant $\tl d(g')=d(x')=1$ dans les cas non mentionn\'es dans ce paragraphe.

\subsection{} L'\'etape suivante consiste \`a remplacer $N_{E_{\vert L}}({\goth o}_{E}^\x)\tl J'$ par $\tl J'$. Elle ne concerne que le cas (r-r). 

\begin{Lemme} Dans le cas (r-r), 
\begin{equation}\label{sansnom}
\sum_{g'\in {\scr C}_{\tau,N_{E_{\vert L}}({\goth o}_{E}^\x)\tl J'}^{st}(g)}\frac{d}{\tl d(g')} \tl c(g')\tr \tl \Lambda_{st}(g'\tau)=
\sum_{g'\in {\scr C}_{\tau,\tl J'}^{st}(g)}\frac{d}{\tl d(g')} \tl c(g')\tr \tl \Lambda_{st}(g'\tau).
\end{equation}
\end{Lemme}

\begin{proof} Consid\'erons l'application $\Psi$ de ${\scr C}_{\tau,\tl J'}^{st}(g)$ dans ${\scr C}_{\tau,N_{E_{\vert L}}({\goth o}_{E}^\x)\tl J'}^{st}(g)$ qui \`a la classe de $\tau$-$\tl J'$-conjugaison de $g'\in {\scr Cl}_{\tau}^{st}(g)\cap \tl J'$ associe la classe de $\tau$-$N_{E_{\vert L}}({\goth o}_{E}^\x)\tl J'$-conjugaison de $g'$. Etudions d'abord son image.\\
Soit $g'\in N_{E_{\vert L}}({\goth o}_{E}^\x)\tl J'$. L'ensemble des classes \`a droite modulo $\tl J'$ dans  $N_{E_{\vert L}}({\goth o}_{E}^\x)\tl J'$ est repr\'esent\'e par $\{ \zeta^{2r}, r\in \[ 0,\frac{q-1}{2}\]\}$. Il existe donc $r\in \[ 0,\frac{q-1}{2}\]$ et $h\in \tl J'$ tels que : $g'=\zeta^{2r} h=\zeta^r(\zeta^r h\zeta^{-r})\tau(\zeta^{-r})$.  L'\'el\'ement $\zeta^r h\zeta^{-r}$ appartient \`a $\tl J'$ tandis que $\zeta^r$ appartient \`a $N_{E_{\vert L}}({\goth o}_{E}^\x)\tl J'$ si et seulement si $r$ est pair. Donc si $r$ est pair, $g'$ appartient \`a $\im \Psi$. R\'eciproquement, les classes de $\tau$-$\tl J'$-conjugaison dans la classe de $\tau$-$N_{E_{\vert L}}({\goth o}_{E}^\x)\tl J'$-conjugaison de $g'$ ont un repr\'esentant de la forme $\zeta^{-2s}g'\tau(\zeta^{2s})=\zeta^{2r-4s}(\zeta^{2s}h\zeta^{-2s})$ pour un $s\in \[ 0,\frac{q-1}{2}\]$. L'une d'entre elles est contenue dans $\tl J'$ si et seulement il existe $s\in \[ 0,\frac{q-1}{2}\]$ tel que $\zeta^{2r-4s}\in F_{0}^\x$, c'est-\`a-dire $r-2s\equiv 0 \mod (\frac{q+1}{2})$.\\
Si $q\equiv 1 \mod 4$, $g'$ est donc $\tau$-$N_{E_{\vert L}}({\goth o}_{E}^\x)\tl J'$-conjugu\'e \`a un \'el\'ement de $\tl J'$ : $\Psi$ est surjective.\\
Si $q\equiv -1 \mod 4$, $g'$ est  $\tau$-$N_{E_{\vert L}}({\goth o}_{E}^\x)\tl J'$-conjugu\'e \`a un \'el\'ement de $\tl J'$ si et seulement si $r$ est pair. \\
Etudions les fibres de $\Psi$. Soient $g_{1}, g_{2}$ deux \'el\'ements de $\tl J'$ qui sont $\tau$-$N_{E_{\vert L}}({\goth o}_{E}^\x)\tl J'$-conjugu\'es. Quitte \`a $\tau$-$\tl J'$-conjugu\'e $g_{1}$, on peut supposer qu'il existe $r\in\[ 0,\frac{q-1}{2}\]$ tel que $\zeta^{2r}$ $\tau$-conjugue $g_{2}$ en $g_{1}$. Mais alors $\zeta^{-4r}$ appartient \`a $\tl J'$, c'est-\`a-dire $\frac{q+1}{2}$ divise $2r$.\\
Si $q\equiv 1 \mod 4$, $\frac{q+1}{2}$ divise $r$ et $g_{1}=g_{2}$ : $\Psi$ est injective. En remarquant que $\Psi$  conserve la valeur de $\frac{d}{\tl d(g')}\tl c(g')\tr \tl \Lambda_{st}(g'\tau)$, on \'etablit l'\'egalit\'e (\ref{sansnom}).\\
Si $q\equiv -1 \mod 4$, $\frac{q+1}{4}$ divise $r$ et $g_{2}$ est \'egal \`a $g_{1}$ ou $\zeta^{-\frac{q+1}{2}}g_{1}\tau(\zeta^{\frac{q+1}{2}})$ : les fibres de $\Psi$ sont de cardinal 2. \\
Dans ce cas, on remarque que l'application $g' \mapsto \zeta^{-1}g'\tau(\zeta)$ d\'efinit une bijection de $\im \Psi$ sur son compl\'ementaire dans ${\scr C}_{\tau,N_{E_{\vert L}}({\goth o}_{E}^\x)\tl J'}^{st}(g)$ qui pr\'eserve la valeur de $\frac{d}{\tl d(g')}\tl c(g')\tr \tl \Lambda_{st}(g'\tau)$. On retrouve donc l'\'egalit\'e (\ref{sansnom}).
\end{proof}

\subsection{} Il reste \`a \'etudier la restriction de l'application $\N$ \`a ${\scr C}_{\tau,\tl J'}^{st}(g)$.

\begin{Lemme}\label{cle}
Soient $g$ et $x$ comme pr\'ec\'edemment. \\
(i) Soit $g'\in \tl J'$. Il existe $h\in \tl J'$ tel que $\N(h^{-1}g'\tau(h))\in J'$.\\ 
(ii) Soient $x_{1}, x_{2}$ deux \'el\'ements de $J'$ conjugu\'es sous $\tl J'$. Alors ils sont conjugu\'es sous $J'$.\\
(iii) Soit $x'\in {\scr Cl}^{st}(x)\cap J'$. Il existe $g'\in {\scr Cl}_{\tau}^{st}(g)\cap \tl J'$ tel que : $\N (g')=x'$. Il en existe deux \`a $\tau$-$\tl J'$-conjugaison pr\`es dans le cas 1 et $4$ dans le cas 2.
\end{Lemme}

\begin{proof} On se place dans le cas 1.\\
(i) On a $g'=zj$ o\`u $z\in F^\x$, $j\in \tl J_{1}$. Par le lemme 3.1 de \cite{Bl2}, il existe $h\in \tl J_{1}$ tel que $\N(h^{-1}j\tau(h))\in J_{1}$. Alors $\N(h^{-1}g'\tau(h))\in J'$.\\
(ii) Notons pour $i=1$ ou 2, $x_{i}=u_{i}j_{i}$ o\`u $u_{i}\in F_{\vert F_{0}}^1$, $j_{i}\in J_{1}$. Soit $h\in \tl J'$ tel que $x_{2}=h^{-1}x_{1}h$. On peut prendre $h\in \tl J_{1}$. Alors :
\begin{equation*}
x_{2}=h^{-1}x_{1}h \ssi u_{1}^{-1}u_{2}=h^{-1}j_{1}hj_{2}^{-1} \in F_{\vert F_{0}}^1\cap \tl J_{1}\ssi u_{1}=u_{2} \text{ et } j_{2}=h^{-1}j_{1}h.
\end{equation*}
Par le lemme 3.1(ii) d\'ej\`a cit\'e, on peut choisir $h\in J_{1}$.\\
(iii) Notons $x'=uj$ avec $u\in F_{\vert F_{0}}^1$ et $j\in J_{1}$. Par le th\'eor\`eme de Hilbert 90 et encore le lemme 3.1(ii), il existe $g'=z\tl j$, $z\in F^\x$ et $\tl j\in \tl J_{1}$ tel que : $\N(g')=x'$. De plus, $g'\in {\scr Cl}_{\tau}^{st}(g)$.\\
Soient maintenant $g'_{1}, g'_{2}\in {\scr Cl}_{\tau}^{st}(g)\cap \tl J'$ tels que : $\N(g_{1})=\N(g_{2})=x'$. Alors $g_{i}$, $i=1,2$, s'\'ecrit $z_{i}j_{i}$ o\`u $z_{i}\in F^\x$ et $j_{i}\in \tl J_{1}$ tel que $j_{i}\tau(j_{i})\in J_{1}$. L'\'egalit\'e pr\'ec\'edente donne alors :
\begin{equation*}
\frac{z_{1}\ol z_{2}}{\ol z_{1}z_{2}}=j_{2}\tau(j_{2})(j_{1}\tau(j_{1}))^{-1} \ssi z_{1}\ol z_{2}\in F_{0}^\x \text{ et } j_{1}\tau(j_{1})=j_{2}\tau(j_{2})
\end{equation*}
quitte \`a multiplier $z_{1}$ et $j_{1}$ par un \'el\'ement de $1+{\goth p}$ sans changer leur produit.
Encore une fois, $j_{1}$ et $j_{2}$ sont donc $\tau$-$\tl J_{1}$-conjugu\'es. Alors $g_{1}$ et $g_{2}$ sont $\tau$-$\tl J'$-conjugu\'es si et seulement si $z_{1}\ol z_{2}$ est une norme de $F^\x$ dans $F_{0}$. A $\tau$-$\tl J'$-conjugaison pr\`es, on a deux $g'\in \tl J'$ de norme cyclique $x'$. 

On se place dans le cas 2. Les assertions (i) et (ii) sont imm\'ediates. Soit $x'\in {\scr Cl}^{st}(x)\cap J'$. Alors il existe $j\in \tl J_{1}$,  unique \`a $\tau$-$\tl J_{1}$-conjugaison pr\`es, tel que $\N (j)=x'$.  Mais pour tout $u\in F_{0}^\x$, $uj$ est de norme cyclique $x'$ et $uj$ est $\tau$-$\tl J'$-conjugu\'e \`a $j$ si et seulement si $u$ est un carr\'e dans $F_{0}^\x$.\\
En effet, si $u$ est un carr\'e, il est clair que $j$ et $uj$ sont $\tau$-$\tl J'$-conjugu\'es. R\'eciproquement,  on peut supposer que $u$ est 1, $\varpi_{0}, \zeta^{q+1}$ ou $\varpi_{0}\zeta^{q+1}$. Si  $j$ et $uj$ sont $\tau$-$\tl J'$-conjugu\'es, il existe $v\in F_{0}^\x$ et $j'\in \tl J_{1}$ tels que $uj=v^{-1}({j'}^{-1}j\tau(j'))\tau(v)$. Mais alors, $uv^2\in \tl J_{1}\cap F_{0}^\x=1+{\goth p}_{0}$ donc $u=1$.
\end{proof}

\subsection{}  Etablissons l'\'egalit\'e  (\ref{id}). Si ${\scr Cl}_{\tau}^{st}(g)\cap \tl J'$ est vide, ${\scr Cl}^{st}(x)$ ne rencontre pas $J'$ (lemme \ref{cle}, (i)) donc les deux sommes de (\ref{tr2}) sont nulles. Supposons maintenant que $g\in \tl J'$ et $x\in J'$ (lemme \ref{cle}, (i)). Alors l'application $\Nc$ induit une surjection de ${\scr C}_{\tau,\tl J'}^{st}(g)$ sur ${\scr C}_{J'}^{st}(x)$ dont les fibres sont de cardinal 2 dans le cas 1, 4 dans le cas 2 (lemme \ref{cle}). De plus, $\tr \tl \Lambda$ est constant sur ces fibres et l'identit\'e s'obtient comme dans \cite[\S 4.5 et 4.6]{Bl2} en ajoutant le r\'esultat suivant :

\begin{Lemme} Soit $g'\in \tl J'$ tel que $x':=\N (g') \in J'$. Alors 
$\tl d(g')$ est \'egal \`a $d(x')$ dans le cas 1, $2d(x')$ dans le cas 2.
\end{Lemme} 

\begin{proof}
Le r\'esultat est imm\'ediat dans les cas (nr-nr) et (r-r) quand $q\equiv 1 \mod 4$.\\
 Dans le cas (r-r) quand $q\equiv -1 \mod 4$, $J$ est le produit $\{\pm 1\} J'$ donc $d(x')$ vaut 1. D'autre part, $N_{E_{\vert L}}({\goth o}_{E}^\x)\tl J'$ est un sous-groupe d'indice 2 de $\tl J^+$ donc la classe de $\tau$-$\tl J^+$-conjugaison de $g'$ contient au plus deux classes de $\tau$-$N_{E_{\vert L}}({\goth o}_{E}^\x)\tl J'$-conjugaison, \`a savoir celles de $g'$ et de $-\zeta^{-1}g'\tau(\zeta)$. Remarquons qu'il existe un entier $r$ tel que $g'\in \zeta^{(q+1)r}\varpi_{0}^{\Bbb Z}\tl J_{1}$ et que toute la classe de $\tau$-$N_{E_{\vert L}}({\goth o}_{E}^\x)\tl J'$-conjugaison de $g'$ est contenue dans $\displaystyle \bigcup_{s\equiv 0  [4]}\zeta^s\varpi_{0}^{\Bbb Z}\tl J_{1}$. Mais,  $-\zeta^{-1}g'\tau(\zeta)$ appartient \`a $\zeta^s\varpi_{0}^{\Bbb Z}\tl J_{1}$ o\`u $s=\frac{q^2-1}{2}-2+(q+1)r\equiv 2 \mod 4$. Par suite : $\tl d(g')=2=2d(x')$.\\
 Dans le cas (r-nr), $g'$ s'\'ecrit $zg'_{1}$ avec $z\in (F^\x)^+$ et $g'_{1} \in \tl J_{1}$. Alors $x'=\frac{z}{\ol z}x'_{1}$ o\`u $x'_{1}=\N (g'_{1})\in J_{1}$. Puisque $z$ et $\frac{z}{\ol z}$ sont centraux dans $\tl J^+$ et $J$ respectivement, on a : $\tl d(g')=\tl d(g'_{1})$ et $d(x')=d(x'_{1})$. On peut donc supposer que $g'\in \tl J_{1}$ et  $x'\in J_{1}$.  Sous cette hypoth\`ese, $\tl J'$ contient toute la classe de $\tau$-$\tl J^+$-conjugaison de $g'$. D'apr\`es le lemme \ref{cle}, $\N$ induit une application de l'ensemble des classes de $\tau$-$\tl J'$-conjugaison contenues dans ${\scr Cl}_{\tau-\tl J^+}(g')$ dans l'ensemble des classes de $J'$-conjugaison contenues dans ${\scr Cl}_{J}(x')$. Cette application est clairement surjective car tout $x''\in {\scr Cl}_{J}(x')$ s'\'ecrit $j^{-1}x'j$, $j\in J$, et est l'image par $\N$ de $g'':=j^{-1}g'\tau(j)$, \'el\'ement appartenant \`a ${\scr Cl}_{\tau-\tl J^+}(g')\cap  \tl J_{1}$. \\
 Etudions la fibre en $x''$. D'apr\`es le lemme \ref{cle}, elle contient au plus 2 classes si $q\equiv 1  \mod 4$, 4 si $q\equiv -1  \mod 4$ et, d'apr\`es la d\'emonstration de ce m\^eme lemme, ces classes de $\tau$-$\tl J'$-conjugaison sont repr\'esent\'ees par : $g''$ et $\alpha_{0}g''$ si   $q\equiv 1  \mod 4$~; $ug''$ o\`u $u$ parcourt un ensemble de repr\'esentants des classes de $F_{0}^\x$ modulo ses carr\'es. Mais, une de ces classes est contenue dans ${\scr Cl}_{\tau-\tl J^+}(g')$  si et seulement si son repr\'esentant $zg''$ est $\tau$-$\tl J^+$-conjugu\'e \`a $g''$, c'est-\`a-dire il existe $h=(\varpi \zeta)^r \zeta^{2s}y$, $r,s\in \Bbb Z$ et $y\in \tl J_{1}$ tel que :
 \begin{equation*}\begin{aligned}
&zg''=((\varpi \zeta)^r \zeta^{2s}y)^{-1} g'' \tau((\varpi \zeta)^r \zeta^{2s}y) \ssi z=\varpi_{0}^{-r}\zeta^{-(q+1)(r+2s)}\mod \tl J_{1}\\ \ssi &z\!=\varpi_{0}^{-r}\zeta^{-(q+1)(r+2s)}\!\! \!\mod 1\!+\!{\goth p}_{0}\! \ssi 
\!\begin{cases}
 r=0 \text{ et } z=1\!&\!\text{si }\! q\equiv 1\!\mod 4\\
 z=1\! \text{ ou }\! z=\varpi_{0}\zeta^{q+1}\!\! &\!\!\text{si }\! q\!\equiv \!-1\!\!\mod 4\\
\end{cases}
\end{aligned}\end{equation*}
Ainsi, la fibre en $x''$ est de cardinal 1 si $q\equiv 1 \mod 4$, 2 sinon. 
\end{proof}

A ce point, sont d\'emontr\'ees les assertions du th\'eor\`eme \ref{BCdim2} correspondant aux cas (nr-nr) et (r-r) lorsque $q$ est congru \`a 1 modulo 4.

\subsection{} Dans les autres cas de niveau strictement positif, chaque repr\'esentation $\tl \pi_{r}$ est l'image par le changement de base stable d'une repr\'esen\-ta\-tion $\pi_{r}''$ \cite[\S 11.4.1]{Ro} caract\'eris\'ee par l'identit\'e (\ref{id}). D'apr\`es ce qui pr\'ec\`ede, pour tout $x\in G$ r\'egulier, elliptique, de la forme $\N (g)$ pour un $g\in \tl G$,
\begin{equation*}
\sum _{i=0}^{d-1} \tr \pi ''_i(x)=\sum _{i=0}^{d-1}\tr \pi_i(x).
\end{equation*}
Alors, comme au paragraphe 4.7 de \cite{Bl2}, on obtient que les images des repr\'esen\-ta\-tions $\pi_{r},  0\leq r\leq d-1$, par le changement de base stable sont les repr\'esentations $\tl \pi_{r}$, $0\leq r\leq d-1$. \\
Dans le cas (r-r) o\`u $q$ est congru \`a -1 modulo 4, les deux repr\'esentations $\tl \pi_{0}$ et $\tl \pi_{1}$ ne sont autres que $\tl\pi\cdot \mu^{-1}\circ \Det $ et $\tl \pi'\cdot \mu^{-1}\circ \Det $ et seule $\tl \pi'\cdot \mu^{-1}\circ \Det $ a pour caract\`ere central $\tl \omega_{\pi}$. Elle est donc l'image de $\pi$ par le changement de base stable.\\
Dans le cas (r-nr), lorsque $q$ est congru \`a -1 modulo 4, l'argument du ca\-rac\-t\`ere central permet d'affirmer que l'image de $\pi$ par le changement de base stable appartient \`a $\{ \tl \pi_{r}, r\equiv 0 \mod 2\}$. Sans hypoth\`ese sur $q$, pour montrer que l'image de $\pi$ est bien $\tl \pi_{0}$, il suffit d'\'evaluer leurs caract\`eres en des \'el\'ements bien choisis. Comme dans  la d\'emonstration de la proposition \ref{stable}, on choisit $g=\zeta$ ou $\zeta^2$. Alors, par (\ref{tr-tlpi(iii)}) et en poursuivant le raisonnement comme en \cite[\S 4.7]{Bl2} :
 \begin{equation*}
\tr \tl \pi(g\tau)=\tr \lambda(x)+\tr \lambda(x^{-1}) =\tr \pi(x)\neq 0.
\end{equation*}
Ceci termine la d\'emonstration du th\'eor\`eme dans le cas des repr\'esentations de niveau strictement positif.
 
\subsection{}  Consid\'erons le cas o\`u $\pi$ est de niveau 0. La repr\'esentation $\tl \pi$ est une repr\'esentation $\tau$-invariante, \`a caract\`ere central trivial sur $F_{0}$ et de caract\`ere stable : elle appartient donc \`a l'image du changement de base stable et son ant\'ec\'edent est n\'ecessairement une repr\'esentation de niveau 0 d\'ecrite par \ref{donnees2} (a),
cons\'equence de tout ce qui pr\'ec\`ede. Ces derni\`eres se distinguent par la restriction de leurs caract\`eres en les \'el\'ements de $G$ dont le groupe des points sur $F$ du centralisateur est le groupe multiplicatif d'une extension quadratique de $F$, \'el\'ements dits ``elliptiques'' dans \cite{Fr} et ``tr\`es elliptiques'' ici pour \'eviter les confusions. 

Soit $g\in \tl G$ tel que $x=g\tau(g)$ est un \'el\'ement de $G$ r\'egulier et tr\`es elliptique. Le groupe $T(F)$ des points sur $F$ du centralisateur de $x$ est le groupe multiplicatif d'une extension quadratique $E$ de $F$, n\'ecessairement non ramifi\'ee sur $F$.

Si ${\scr Cl}_{\tau}^{st}(g)\cap \tl J$ est vide, alors ${\scr Cl}^{st}(x)\cap J$ aussi. En effet, un \'el\'ement $x'$ de ${\scr Cl}^{st}(x)\cap J$ est de la forme $\N (g')$ o\`u $g'\in {\scr Cl}_{\tau}^{st}(g)$ et commute \`a $x'$ (puisque $g$ commute \`a $x$). Or $x'$ est tr\`es elliptique donc $g'$ appartient \`a $F^\x\tl U_{0}(\scr L)=\tl J$ \cite[th. 1]{Fr}. Pour un tel $g$, l'identit\'e (\ref{id}) est satisfaite.

Si ${\scr Cl}_{\tau}^{st}(g)\cap \tl J$ n'est pas vide, on consid\`ere $g'\in {\scr Cl}_{\tau}^{st}(g)\cap \tl J$ et $x'=\N (g')$. Alors $x'\in \tl U_{0}(\scr L)$ et $x'$ est $\tl G$-conjugu\'e \`a $x$ : $x'=y^{-1}xy, y\in \tl G$.\\
Notons que $\tau(x')$ appartient \`a $\tl U_{0}(\scr L)$ et $\tau(x')=(y^{-1}\tau(y))^{-1}x'y^{-1}\tau(y)$. Ainsi $y^{-1}\tau(y)\in F^\x\tl U_{0}(\scr L)$, autrement dit : $\tau(y)=yj$, $j\in \tl J$. Mais alors la cha\^\i ne de r\'eseaux $y\scr L$ est une cha\^\i ne autoduale de $F^2$ et $x\in U_{0}(y\scr L)$. De deux choses l'une :\\
- la cha\^\i ne $y\scr L$ est de m\^eme invariant que $\scr L$ auquel cas il existe $y_{0}\in G$ tel que $y\scr L=y_{0}\scr L$~: $y_{0}^{-1}xy_{0}\in U_{0}(\scr L)=J$ ;\\
- la cha\^\i ne $y\scr L$ n'est pas de m\^eme invariant que $\scr L$ auquel cas  il existe $y_{0}\in G$ tel que $y\scr L=y_{0}{\scr L}'$ o\`u ${\scr L}'$ est une cha\^\i ne de p\'eriode 1  telle que $U_{0}({\scr L}')$ contient le m\^eme sous-groupe d'Iwahori ${\scr I}$ que $U_{0}(\scr L)$ : $y_{0}^{-1}xy_{0}\in U_{0}({\scr L}')$. Mais, ${\scr L}'$ est  d'invariant pair et $U_{0}({\scr L}')$ est la r\'eunion disjointe de ${\scr I}$ et de son compl\'ementaire, le premier form\'e des \'el\'ements de d\'eterminant dans $1+{\goth p}$, le second form\'e des \'el\'ements de d\'eterminant dans $-1+{\goth p}$. Or le d\'eterminant de $y_{0}^{-1}xy_{0}$ est \'egal \`a $\Det g'/\ol{\Det g'}$ o\`u 
$\Det g'$ est n\'ecessairement de valuation paire, donc est un \'el\'ement de $1+{\goth p}$. Ainsi, 
$y_{0}^{-1}xy_{0}$ appartient \`a ${\scr I}\subset U_{0}(\scr L)$.\\
Dans les deux cas, ${\scr Cl}^{st}(x)\cap J$ n'est pas vide. On peut donc choisir $x\in J$ puis $g\in \tl J$. Par la formule de Mackey et puisque $x$ est tr\`es elliptique, on obtient : 
 \quad  $\tr \tl \pi (g\tau)=\tr \tl \Lambda (g\tau)\quad$ et $\quad \tr\pi(x)=\tr \lambda(x)$.\\
On conclut gr\^ace \`a  (\ref{tr-lambda-0}).

\section{Caract\`eres de repr\'esentations irr\'eductibles de certaines extensions de groupes finis.}\label{Traces}

On pr\'esente dans ce paragraphe un calcul de caract\`eres de repr\'esentations admissibles de groupes compacts qui permet de justifier ou de se convaincre de la validit\'e d'affirmations contenues dans le paragraphe \ref{donnees} (cas (b')) et la d\'emonstration du lemme \ref{construction-choix}.

Le calcul expos\'e est plus complexe que n\'ecessaire ici. Mais, on rencontre cette m\^eme question dans l'\'etude du transfert de $U(1,1)(F_{0})\x U(1)(F_{0})$ \`a $U(2,1)(F_{0})$. On a donc choisi un cadre suffisament grand pour r\'eutiliser les r\'esultats lors de cette \'etude.

\subsection{Donn\'ees, hypoth\`eses et objectif.}\label{Tr1}
On fixe un nombre premier $p$.
 
Soient ${\Bbb J}$ un groupe fini, 
$\J_{1}$ un sous-groupe distingu\'e de ${\Bbb J}$ de centre $Z$ et $\J '$ un sous-groupe de $\J$ tels que $\J=\J' \J_{1}$ et $\J'\cap \J_{1}\subset Z$. On suppose que :

\begin{itemize} 
\item[(i)] $\J_{1}$ est un $p$-groupe, extra-sp\'ecial de classe 2 ou ab\'elien ;
\item[(ii)] $\J'=T\J'_{1}$ o\`u $T$ un sous-goupe ab\'elien de ${\Bbb J}'$ et  $\J'_{1}$ un $p$-sous-groupe distingu\'e, ou bien extra-sp\'ecial de classe 2 dont le centre $Z'$ contient $\J'_{1}\cap T$,  ou bien ab\'elien (dans ce cas on le note aussi  $Z'$)~;
\item[(iii)] $[\J'_{1}, \J_{1}]=1$ ;
\item[(iv)] $\J'/\J'_{1}$ est d'ordre premier \`a $p$ ;
\item[(v)] $Z'Z$ est contenu dans le centre de $\J$ et $\J'_{1}\J_{1}$ est un $p$-groupe extra-sp\'ecial de classe 2, de centre $Z'Z$.
\end{itemize}

Soient  $\theta$ un caract\`ere fid\`ele de $Z$ et $\theta'$ un caract\`ere fid\`ele de $Z'$ tels que $\theta_{\vert Z\cap Z'}=\theta'_{\vert Z\cap Z'}$.   

Ainsi, le groupe $Z$ est un groupe cyclique d'ordre $p$ et le quotient $V=\J_{1}/Z$ est un ${\Bbb F}_{p}$-espace vectoriel de dimension paire, muni d'une forme bilin\'eaire altern\'ee non d\'eg\'en\'er\'ee $<.,.>$ d\'efinie par :
$$\begin{matrix} <.,.>&:&V\x V& \fl &\mu_{p}\\ & & (gZ,g'Z)& \mapsto &\theta([g,g']).
\end{matrix}$$
Ceci est \'egalement valable pour $Z'$ ``en primant'' $V$, $\J_{1}$, $Z$ et $\theta$.

\smallskip	

\noindent  On note $\Theta$ le caract\`ere $\theta'\cdot \theta$ de $Z'Z$ et $\eta_{\Theta}$ la repr\'esentation de Heisenberg de $\J'_{1}\J_{1}$ de caract\`ere central $\Theta$. Elle est isomorphe au produit tensoriel des repr\'esentations de Heisenberg de $\J_{1}$ et $\J'_{1}$ de caract\`eres centraux $\theta$ et $\theta'$ respectivement. On note $\eta_{\theta}$ et $\eta_{\theta'}$ ces deux repr\'esentations, $p^a$ et $p^b$ leurs dimensions respectives.

\smallskip

\noindent On suppose que les repr\'esentations $\eta_{\Theta}$ et $\eta_{\theta'}$ se prolongent en des repr\'esen\-tations de $\J$ et $\J'$ respectivement.

\medskip	

Notre but est de d\'efinir une bijection entre les prolongements de $\eta_{\Theta}$ \`a $\J$ et ceux de $\eta_{\theta'}$ \`a $\J'$ caract\'eris\'ee par une identit\'e de caract\`eres. 

\noindent Dans les deux cas, les prolongements se distinguent par les valeurs de leurs caract\`eres sur les \'el\'ements de $T$. On cherche donc \`a calculer la trace des prolongements de $\eta_{\Theta}$ en un \'el\'ement de $T$ en fonction de celle des prolongements de $\eta_{\theta'}$, au moins dans les situations qui nous sont utiles.

\subsection{Cinq cons\'equences des donn\'ees et hypoth\`eses.}\label{Tr2}
\subsubsection{ } \label{Tr21}Tout d'abord, de l'hypoth\`ese (iv), on d\'eduit que :
$$\J'\cap \J_{1}\subset \J'_{1}\quad \text{ d'o\`u } \quad \J/\J'_{1}\J_{1}\simeq \J'/\J'_{1}\simeq T/T\cap Z'\quad \text{ et } \quad \J'\cap \J_{1}\subset Z',$$
le dernier isomorphisme et l'inclusion sont alors cons\'equences de (ii) et (i) respectivement.

\subsubsection{ } \label{Tr22} A $x\in \J'$, on associe l'automorphisme symplectique de $V$, encore not\'e $x$, d\'efini par la conjugaison par $x$, le sous-espace $V^x$  de $V$ form\'e des points fixes sous $x$ et le caract\`ere $\chi_{x}$ de $V^x$ d\'efini par : 
$$\chi_{x}(v)=\theta([x,v]) \text{ pour tout } v\in V^x.$$
Remarquons d'une part, que $\J'_{1}$ centralisant $\J_{1}$ (iii), $V^x$ et $\chi_{x}$ ne d\'ependent que de la classe de $x$ modulo $\J'_{1}$~; d'autre part, que si $x$ n'appartient pas \`a $\J'_{1}$, son ordre $r$  dans $\J'/\J'_{1}$ est premier \`a $p$ (par (iv)) et $\chi_{x}$ est trivial. Pour cette derni\`ere assertion, il suffit de remarquer que :
\begin{itemize}
\item[-] pour tout $v\in V^x$, tout $s\in {\Bbb N}$, $\chi_{x}(v)^s=1 \ssi s\equiv 0 [p]\, \text{ ou } \chi_{x}(v)=1$; 
\item[-] pour tout $v\in V^x$, $\chi_{x}(v)^r=\chi_{x^r}(v)=1$.
\end{itemize} 

\noindent Par cons\'equent, puisque $\theta$ est fid\`ele, la projection naturelle du sous-groupe $\J_{1}^x$ des \'el\'ements de $\J_{1}$ invariants par conjugaison par $x$ sur $V^x$ est surjective. En primant les notations, on obtient le r\'esultat analogue pour tout \'el\'ement de $T$. Ainsi,
\begin{equation}\label{(Tr221)}
\forall x\in \J', \quad \left( \J_{1}/Z\right)^x\simeq\J_{1}^x/Z \quad \text{ et } \quad \forall t\in T, \quad \left( \J'_{1}/Z'\right)^t\simeq{\J'_{1}}^t/Z'\end{equation}
Il s'en suit que : 
\begin{equation*}
\forall t\in T, \quad \left( \J'_{1}\J_{1}/Z'Z\right)^t\simeq{\J'_{1}}^t\J_{1}^t/Z'Z.\end{equation*}

\subsubsection{ } \label{Tr23} Soit $t\in T$. L'application de $V$ dans $V$ qui \`a $v\in V$ associe $[t^{-1},v]$ est lin\'eaire de noyau $V^t$. Son image ${\scr I}_{t}=\{ [t^{-1}, v], v\in V\}$ est l'orthogonal de $V^t$.
De plus, l'action de $T$ par conjugaison sur $\J_{1}$ fournit une repr\'esentation de $T$ dans $GL(V)$ triviale sur $T\cap \J'_{1}$. Puisque $T/T\cap \J'_{1}$ est fini d'ordre premier \`a $p$ (iv), cette repr\'esentation est semi-simple donc somme de caract\`eres. On en d\'eduit que les sous-espaces $V^t$ et ${\scr I}_{t}$ sont suppl\'ementaires. Par cons\'equent,  pour tout $t\in T$ (et par suite pour tout $x\in \J'$), la restriction de $<.,.>$ au sous-espace $V^t$ est non d\'eg\'en\'er\'ee et la dimension de $V^t$ est paire, not\'ee $2a_{t}$. 

\subsubsection{ }\label{Tr24} Les paires $(\J, \lambda)$ et $(\J', \lambda')$ sont deux cas de la situation A.1.7 de \cite{BH2}. En utilisant le corollaire A.1.8 et la surjectivit\'e des projections de $(\J'_{1}\J_{1})^t$ sur $\left( \J'_{1}\J_{1}/Z'Z\right)^t$ et de ${\J'_{1}}^t$ sur $\left( \J'_{1}/Z'\right)^t$ pour tout $t\in T$ (comme dans la d\'emonstration de la proposition 13.1 de \cite{BH1}), on obtient :

\begin{Lemme} Soit $t\in T$.
\begin{itemize} 
\item[a)] $\tr \lambda(t)\not =0, \quad \tr \lambda'(t)\not =0 \quad$ et $\quad \Vert \tr \lambda(t)\Vert=\vert V^t\vert^{1\over 2}\cdot \vert {V'}^t\vert^{1\over 2},\break \quad \Vert \tr \lambda'(t)\Vert=  \vert {V'}^t\vert^{1\over 2}.$
\item[b)] Soit $y\in \J'_{1}\J_{1}$. Alors $\tr \lambda(ty)$ est non nul si et seulement si $ty$ est $\J'_{1}\J_{1}$-conjugu\'e \`a un \'el\'ement de $tZ'Z$.
\item[c)] Soit $y'\in \J'_{1}$. Alors $\tr \lambda'(ty')$ est non nul si et seulement si $ty'$ est $\J'_{1}$-conjugu\'e \`a un \'el\'ement de $tZ'$.
\end{itemize}
\end{Lemme}	

\subsubsection{ }  Enon\c cons une derni\`ere cons\'equence.
\begin{Lemme} \label{Tr25} Deux \'el\'ements de $\J'$ sont $\J$-conjugu\'es si et seulement s'ils sont $\J'$-conjugu\'es.\end{Lemme}

\begin{proof} Soient $x, x'$ deux \'el\'ements de $\J'$ qui sont $\J$-conjugu\'es. Quitte \`a conjuguer l'un d'entre eux par un \'el\'ement de $\J'$, on peut supposer que $x$ et $x'$ sont $\J_{1}$-conjugu\'es. On \'ecrit : $x=ty$ o\`u $t\in T$ et $y\in \J'_{1}$ et $x'=jxj^{-1}$ avec $j\in \J_{1}$.  Montrons que $x'$ et $x$ sont \'egaux. \\
Si $x\in \J'_{1}$, il suffit d'utiliser (iii). Si $x\not \in \J'_{1}$, on a gr\^ace \`a (iii) : 

\quad \quad \quad \quad \quad \quad $x'=ty\cdot y^{-1}[t^{-1},j]y\cdot [y^{-1},j]=x\cdot  [t^{-1},j]$.\\
  Ainsi, $[t^{-1},j]=x^{-1}x'$ appartient \`a $Z$ ($\J_{1}\cap \J' \subset Z$ par hypoth\`ese), c'est-\`a-dire que  $j\in V^t=V^x$. Et puisque $\chi_{x}$ est trivial (\S \ref{Tr22}) et $\theta$ fid\`ele, on conclut que $x$ et $x'$ sont \'egaux.\end{proof}

\noindent Par  les lemmes \ref{Tr24} et \ref{Tr25}, on obtient :

\begin{Corollaire} Soit $x$ un \'el\'ement de $\J'$. Les quatre propositions suivantes sont \'equivalentes :
\begin{itemize}
\item[\it (i)] $\tr \lambda (x)\not =0$ ;
\item[\it (ii)] $\tr \lambda' (x)\not =0$ ;
\item[\it (iii)] $x$ est $\J'$-conjugu\'e \`a un \'el\'ement de $TZ'$ ;
\item[\it (iv)] $x$ est $\J$-conjugu\'e \`a un \'el\'ement de $TZ'$.
\end{itemize}
\end{Corollaire}

\subsection{Comparaison de traces.}\label{Tr3} 
\subsubsection{ } \label{Tr31} Soit $\lambda$ un prolongement de $\eta_{\Theta}$ \`a $\J$. Notons  $\J'_{0}$ le centralisateur de $\J_{1}$ dans $\J'$. La restriction de $\lambda$ \`a $\J'_{0}$ est de la forme $\lambda'_{0}\cdot \eta_{\theta}$ o\`u $\lambda'_{0}$ est un prolongement de $\eta_{\theta'}$ \`a $\J'_{0}$. Soit  $\lambda'$ un prolongement  de $\lambda'_{0}$ \`a $\J'$.

\smallskip	

\noindent
Par la suite, on note  ${\scr X}$ l'ensemble des caract\`eres de $\J$ triviaux sur $\J'_{0}\J_{1}$. Il s'identifie via la restriction \`a l'ensemble des caract\`eres de $\J'$ triviaux sur $\J'_{0}$ ou \`a ceux de $T$ triviaux sur $T\cap \J'_{0}$. Selon le contexte, on regarde ${\scr X}$ d'une fa\c con ou d'une autre et on note $d$ son cardinal.

\smallskip	

\noindent Remarquons que lorsque $a=0$, $\lambda$ est isomorphe \`a $\lambda'\otimes \theta$ pour un prolongement $\lambda'$ de $\eta_{\theta'}$ convenable. Il existe donc une bijection  que l'on peut d\'efinir ``canoniquement'' par l'\'egalit\'e des traces sur les \'el\'ements de $T$. 

\noindent Si $d=1$, $\J$ n'est autre que le produit $\J'_{0}\J_{1}$ et $\J'$ est \'egal \`a $\J_{0}'$.   La repr\'esentation $\lambda$ est de la forme  ${\lambda'}\otimes \eta_{\theta}$ pour un unique prolongement $\lambda'$ de $\eta_{\theta'}$. Il existe donc une bijection  que l'on peut d\'efinir ``canoniquement'' par l'identit\'e de caract\`eres suivante :
$$\forall x\in T,\quad \tr \lambda(x)= p^a\tr\lambda'(x).$$
On suppose donc : $a>0$ et $d>1$. 

\subsubsection{ } \label{Tr32} On d\'esigne par $\K$ le groupe $\J'Z$. Puisque $\lambda$ prolonge $\eta_{\Theta}$, sa restriction \`a $\K$ est une sous-repr\'esentation de la restriction \`a $\K$ de $\ind_{\J'_{0}\J_{1}}^\J \lambda'_{0}\otimes \eta_{\theta}$ qui est multiple de $\ind_{\J'_{0}Z}^\K \lambda_{0}'\cdot \theta$ d'apr\`es la formule de Mackey. Il existe donc $d$ entiers  positifs $m_{\xi}$ ($\xi\in {\scr X}$) tels que :
\begin{equation}\label{equ:n321}
\forall k\in \K, \quad \lambda(k)=\bigoplus_{\xi\in {\scr X}} m_{\xi}(\xi \lambda'\cdot \theta)(k)\quad \text{ et } \quad \sum_{\xi\in {\scr X}}m_{\xi}=p^a
\end{equation}
d'o\`u : \quad \quad \quad \quad \quad \quad \quad \quad 
$\displaystyle \forall x\in T, \quad \tr\lambda(x)=\tr\lambda'(x)\cdot \sum_{\xi\in {\scr X}}m_{\xi}\xi(x).$

\smallskip	

\noindent Notre but est d'\'evaluer $\sum_{\xi\in {\scr X}}m_{\xi}\xi(x)$.  De l'assertion {\it a} du lemme \ref{Tr24}, on d\'eduit qu'elle est de module  $\vert V^x\vert^{1/2}$ pour tout $x$ de $T$. On obtient alors un syst\`eme d'\'equations en les inconnues $m_{\xi}$, \`a savoir~: 
\begin{equation}\label{equ:n322}\begin{split}
&\sum_{\xi\in {\scr X}}m_{\xi}=p^a \\ 
&\sum_{\xi\in {\scr X}}m_{\xi}^2 + \sum_{\xi'\in {\scr X}-\{1\}}(\sum_{\xi\in {\scr X}}m_{\xi}m_{\xi\xi'})\xi'(x)=\vert V^x\vert \text{ pour chaque } \ol x\not =1, 
\end{split}\end{equation}
o\`u $\ol x$ repr\'esente la classe de $x\in T$ modulo $\J'_{0}\cap T$.

\subsubsection{ } \label{Tr33} On suppose dor\'enavant que : 

\smallskip	
\noindent {\it (H)}\, {\it $T/\J'_{0}\cap T$ est cyclique et que pour tout $\ol x\in T/\J'_{0}\cap T$ non trivial, $\vert V^x\vert =1$.}
\smallskip	

\noindent Alors les $d$ inconnues $M_{1}=\sum_{\xi\in {\scr X}}m_{\xi}^2 -1$ et $M_{\xi'}=\sum_{\xi\in {\scr X}}m_{\xi}m_{\xi\xi'}$, $\xi'\in {\scr X}$ non trivial, sont solutions du syst\`eme lin\'eaire 
\begin{equation*}
\left\{ \begin{aligned} &M_{1}+ \sum_{\xi'\in {\scr X}-\{1\}}M_{\xi'}=p^{2a}-1\\
& M_{1}+ \sum_{\xi'\in {\scr X}-\{1\}}\xi'(x)M_{\xi'}=0,\\
\end{aligned}\right.\end{equation*}
c'est-\`a-dire : \
$\forall \xi\in {\scr X}, \quad M_{\xi}={p^{2a}-1\over d}$. 
 En particulier : $\displaystyle \sum_{\xi\in {\scr X}}m_{\xi}^2-1={p^{2a}-1\over d}$.

\smallskip	

Soit $\{ m_{\xi}, \xi\in {\scr X}\}$ une solution enti\`ere de (\ref{equ:n322}). Alors :
\begin{equation}\label{equ:n331}
\sum_{\xi\in {\scr X}}\left({p^a+1\over d}-m_{\xi}\right)^2=\vert {\scr X}\vert\left({p^a+1\over d}\right)^2-2{p^a+1\over d}\sum_{\xi\in {\scr X}}m_{\xi}+\sum_{\xi\in {\scr X}}m_{\xi}^2=1
\end{equation}
{\it i)} Si ${p^a+1\over d}\in {\Bbb Z}$, il existe, \`a permutation pr\`es des $m_{\xi}$, deux solutions enti\`eres de (\ref{equ:n331}), \`a savoir : 

\quad \quad \quad \quad \quad \quad \quad \quad \quad \quad 
$m_{1}={p^a+1\over d}\pm 1$ et $m_{\xi}={p^a+1\over d}$ si $\xi\not =1$ ; 

\noindent et, toujours \`a permutation pr\`es des $m_{\xi}$, il existe une unique solution enti\`ere de (\ref{equ:n322}), \`a savoir : 

\quad \quad \quad \quad \quad \quad \quad \quad \quad \quad 
$m_{1}={p^a+1\over d}-1$ et $m_{\xi}={p^a+1\over d}$ si $\xi\not =1$.

\noindent {\it ii)} Si ${p^a+1\over d}\not\in {\Bbb Z}$, chaque terme de la somme (\ref{equ:n331}) appartient \`a $]0,1[$ donc, pour tout $\xi\in {\scr X}$, $m_{\xi}$ est \'egal \`a $u$ ou $u+1$ o\`u $u$ est la partie enti\`ere de ${p^a+1\over d}$. Notons $r$ la diff\'erence $(p^a+1)-ud$ et $n$ le nombre d'entiers $m_{\xi}$ \'egaux \`a $u$. Remarquons que : $0<r<d$.

\noindent De la premi\`ere \'equation de (\ref{equ:n322}), on d\'eduit que $n=d+1-r$, puis de (\ref{equ:n331}), on obtient :
\begin{equation*}\begin{split}
&\sum_{\xi\in {\scr X}}\left({p^a+1\over d}-m_{\xi}\right)^2=(d+1-r)\left({r\over d}\right)^2+(r-1)\left(1-{r\over d}\right)^2=1\\ \ssi & (r-2)(r-d)=0Ê\impl r=2.
\end{split}\end{equation*}
Dans ce cas, $d$ divise donc $p^a-1$, l'un des $m_{\xi}$ est \'egal \`a ${p^a-1\over d}+1$ et tous les autres  \`a ${p^a-1\over d}$. On obtient ainsi toutes les solutions de (\ref{equ:n322}).

\smallskip	

\noindent En cons\'equence :

\begin{Lemme}\label{(H)} On suppose l'hypoth\`ese (H) satisfaite. On note $\J'_{0}$ le centralisateur de $\J_{1}$ dans $\J'$, $d$ le cardinal de $ T/\J'_{0}\cap T$  et $p^{2a}$ celui de $\J_{1}/Z$. 

\noindent Si $d$ divise $p^a+1$, il existe un unique prolongement $\lambda'$ de $\eta_{\theta'}$ tel que :
$$\forall x\in T, \quad \tr \lambda(x)=\begin{cases}
p^a\tr\lambda'(x)&\text{ si } x\in \J'_{0}\cap T,\\
-\tr\lambda'(x)\hfill&\text{ sinon.}\end{cases}$$
\noindent Sinon, il existe un unique prolongement $\lambda'$ de $\eta_{\theta'}$ tel que :
$$\forall x\in T, \quad \tr \lambda(x)=\begin{cases} p^a\tr\lambda'(x)&\text{ si } x\in \J'_{0}\cap T,\\ \tr\lambda'(x)\hfill&\text{ sinon.}\end{cases}$$
\end{Lemme}

\subsubsection{Applications}\label{precisions}\mbox{}\\
(1) L'affirmation de \ref{donnees} (b') demande une justification lorsque la repr\'esenta\-tion $\lambda$ n'est pas de dimension 1. D'apr\`es \cite[\S A.6.2]{Bl1}, $\lambda$ est de la forme $\lambda_{\theta}$ pour un caract\`ere $\theta$ d'un sous-groupe de $J$. On applique ce qui pr\'ec\`ede avec~:
\begin{equation*}
\J=J/\Ker \theta\, , \J_{1}=J\cap U_{1}(\scr L)/\Ker \theta\, , T=H/\Ker \theta\cap H\, , \J'_{1}=Z'=Z.
\end{equation*}
en remarquant que $T/\J'_{0}\cap T\simeq k^1_{\vert k_{0}}$, $V$ est de cardinal $q^2$ et le stabilisateur de chaque \'el\'ement non nul de $V$ est r\'eduit \`a $\{1\}$. 

\smallskip	

\noindent (2) On reprend les notations de \ref{donnees2} (b) et on se place dans le cas (r-nr) du paragraphe \ref{construction}. Pour calculer $\tr \lambda$ sur ${\goth o}_{E_{\vert L}}^1$, on pose :
\begin{equation*}
\J= J/\Ker \theta\, , \J_{1}=J_{1}/\Ker \theta\, , T={\goth o}_{E_{\vert L}}^1/{\goth o}_{E_{\vert L}}^1\cap  \Ker \theta\, \text{ et } \, \J'_{1}=Z'=Z
\end{equation*}
et on applique le lemme pr\'ec\'edent. Remarquons que $\J'_{0}$ est alors \'egal \`a $\{\pm 1\}\J'_{1}$ donc $d$ vaut $\frac{q+1}{2}$ tandis que $a$ vaut 1 ou 0 suivant que la dimension de $\lambda$ est $q$ ou $1$. L'hypoth\`ese (H) est bien satisfaite car $T/T\cap \J'_{0}$ est cyclique et tous les \'el\'ements de ${\goth o}_{E_{\vert L}}^1$ non centraux sont minimaux.

\smallskip	

\noindent (3) On se place dans la m\^eme situation que pr\'ec\'edemment mais on s'int\'eresse cette fois \`a la repr\'esentation $\tl \lambda$. On suppose de plus que $\tl \lambda$ n'est pas de dimension 1 (donc $a=1$). On pose : 
\begin{equation*}
\J= \tl J/\Ker\tl \theta\, , \J_{1}=\tl J_{1}/\Ker \tl \theta\, , T={\goth o}_{E}^\x/{\goth o}_{E}^\x\cap  \Ker\tl \theta\, \text{ et } \, \J'_{1}=Z'=Z.
\end{equation*}
C'est un cas particulier de \ref{Tr1} et toutes les cons\'equences qui suivent sont valables, en particulier \ref{Tr32}.


\end{document}